\documentclass[preprint]{preprint}%review
\usepackage{lineno}

\journal{}

\usepackage{amsmath}
\usepackage{amssymb}
\usepackage{amsthm}
\usepackage{graphicx}
\usepackage{epic,eepic,epsfig}
\usepackage{color}
\usepackage{subfigure}  % add 1.22
\usepackage{placeins}   % add 1.22
\usepackage{multirow}
\usepackage{epstopdf}
\usepackage{varwidth}
\usepackage{float}
\usepackage[table]{xcolor}
\usepackage{subfigure}
\usepackage{pdfpages}

\newtheorem{theorem}{Theorem}[section]
\newtheorem{example}{Example}[section]

\newtheorem{remark}{Remark}[section]

\definecolor{tabclr}{cmyk}{0,0,1,0}
\usepackage{charter}

\usepackage{subfigure}
\usepackage{diagbox}
\usepackage{xurl}
\usepackage{xcolor}
\newcommand{\h}{\boldsymbol}
\begin{document}

\title{Highly efficient norm preserving numerical schemes for micromagnetic energy minimization based on SAV method}

\author[MUST]{Jiayun He}
\ead{3220004597@student.must.edu.mo}
\author[MUST]{Lei Yang\corref{cor}}
\ead{leiyang@must.edu.mo}
\author[XTU]{Jiajun Zhan}
\ead{zhan@xtu.edu.cn}
\cortext[cor]{Corresponding author}
\address[MUST]{School of Computer Science and Engineering, Faculty of Innovation Engineering, Macau University of Science and Technology, Macao SAR 999078, China}
\address[XTU]{School of Mathematics and Computational Science, Xiangtan University, Xiangtan 411105, Hunan, China.}

\begin{abstract}
In this paper, two efficient and magnetization norm preserving numerical schemes based on the scalar auxiliary variable (SAV) method are developed for calculating the ground state in micromagnetic structures.
The first SAV scheme is based on the original SAV method for the gradient flow model, while the second scheme features an updated scalar auxiliary variable to better align with the associated energy.
To address the challenging constraint of pointwise constant magnetization length, an implicit projection method is designed, and verified by both SAV schemes.
Both proposed SAV schemes partially preserve energy dissipation and exhibit exceptional efficiency, requiring two linear systems with constant coefficients to be solved.
The computational efficiency is further enhanced by applying the Discrete Cosine Transform during the solving process.
Numerical experiments demonstrate that our SAV schemes outperform commonly used numerical methods in terms of both efficiency and stability.

\end{abstract}


\maketitle

\section{Introduction}\label{sec1}

The magnetic ground state defines the most stable magnetic order in a material, achieved when the system's Gibbs free energy is minimized. This equilibrium configuration dictates the material's fundamental magnetic properties.
The analysis of magnetic ground states provides a fundamental basis for understanding key morphological features, such as the dimensions and internal structure of domains and their intervening walls \cite{BernadouMDepeyreS02:1018, BernadouMDepeyreS03:86}.
Understanding the magnetic ground state also enables the systematic investigation of how external stimuli---such as temperature, pressure, and applied magnetic fields---govern the equilibrium domain structure \cite{TellezMoraAHeX24:1}. 
These insights are essential for designing more efficient magnetic devices, such as memory and sensor technologies.

To approach the ground state in micromagnetic material, besides numerical simulation of the commonly used dynamic model presented by the Landau-Lifshitz (LL) or Landau-Lifshitz-Gilbert (LLG) equation \cite{YangLChenJR21:110142, LiPYangL23:182}, one effective approach is to solve the Gibbs free energy minimization problem with the local norm constraint. 
A wealth of literature addresses the numerical solution of the Gibbs energy minimization problem.
The steepest descent method, a fundamental gradient-based optimizer, is commonly employed for Gibbs energy minimization. Its performance is highly dependent on the step length selection. For instance, the Barzilai-Borwein method has been introduced as an effective strategy to accelerate convergence \cite{ExlLBanceS14:17D118}. Further applications and developments of the steepest descent method in this specific context are well-documented in the literature \cite{FuruyaAFujisakiJ15:1, OikawaTYokotaH16:056006}.
Among the most prevalent approaches are variants of the conjugate gradient (CG) method. Early work by Bernadou et al. implemented a combination of the CG and finite element methods to investigate domain wall properties, such as thickness and evolution \cite{BernadouMDepeyreS02:1599, BernadouMDepeyreS02:1018, BernadouMDepeyreS03:86, BernadouMDepeyreS04:164}. To overcome the accuracy limitations of standard CG, subsequent research introduced nonlinear CG methods, which incorporate step-length restrictions in the line search \cite{FischbacherJKovacsA17:045310, TanakaTFuruyaA17:7100304}. Furthermore, preconditioned CG methods have been explored to achieve higher computational efficiency \cite{ExlLFischbacherJ19:179, SussDSchreflT00:3282}.
Beyond the aforementioned techniques, the literature encompasses a range of other strategies for addressing the Gibbs energy minimization problem. These include classical numerical optimization approaches, such as the quasi-Newton method \cite{ScholzWFidlerJ03:366} and the augmented Lagrangian method \cite{BernadouMHeS98:512, BernadouMHeS}, as well as mathematical relaxations of the original problem \cite{CarstensenCPraetoriusD05:2633, CarstensenCProhlA01:65}. More recently, machine learning approaches have also emerged as a powerful tool in this domain \cite{SchafferSSchreflT23:170761}.

The scalar auxiliary variable (SAV) method is an efficient numerical approach originally proposed for gradient flow problems \cite{ShenJXuJ18:407}. Its core idea is to introduce a scalar auxiliary variable, transforming the original system into an equivalent form that facilitates the construction of linearly implicit time-stepping schemes. 
A critical feature of these SAV-based schemes is their unconditional energy stability, thereby preserving the energy dissipation structure of the continuous model in the discrete solution.

Owing to its combination of computational efficiency, stability, and implementation flexibility, the SAV method has attracted considerable attention and been extended in several directions. These include rigorous convergence analysis \cite{ShenJXuJ18:2895}, the development of generalized frameworks \cite{ChengQLiuC21:113532}, and other variants \cite{ChengQShenJ18:A3982, JuLLLiX22:66, JiangMSZhangZY22:110954}. The method's robustness and broad applicability are demonstrated by its successful use in diverse physical systems, including fluid dynamics \cite{LiXShenJ20:2465, LiXLShenJ22:141}, phase-field models \cite{AkrivisGLiBY19:A3703, YangXHeX22:114376, ChenCYangX19:35, LiQMeiLQ21:107290}, density functional theory \cite{WangTZhouJ23:719}, and magneto-hydrodynamics \cite{LiXLWangWL22:1026}.

A fundamental tenet of micromagnetics dictates that the magnetization magnitude must remain invariant everywhere below the Curie temperature.
For the dynamic model described by the LL or LLG equation, numerous numerical methods have been designed to satisfy this pointwise length constraint; please refer to \cite{GarciaCervera2007:Review, ZhanJYangL24:1327, HeJYangL24:1179, HeJYang}. 
Specifically, a simple implicit projection method was introduced to solve the LL equation in \cite{EWNWangXP01:1647}.
This scheme has been proven unconditionally stable and temporally first-order convergent.
In approaching the ground state of magnetization, the SAV method is more efficient when considering micromagnetic flow, and the projection method ensures that the magnetization length is preserved everywhere.

In this paper, two efficient numerical schemes are constructed based on the idea of the SAV method to solve the Gibbs energy minimization problem. 
The original SAV method introduces an auxiliary variable for nonlinear terms to avoid expensive computational costs. The significant computation of the demagnetization field, due to its long-range character, prompted us to adopt the same strategy.
Firstly, by introducing a scalar auxiliary variable for the magnetostatic energy and combining the projection method, the first SAV1 scheme is proposed. 
To achieve a more accurate representation of the magnetostatic energy, we introduce the SAV2 scheme with a modified update strategy for the scalar auxiliary variable.
For both SAV schemes, the pointwise constraint of length of magnetization is preserved by projecting the magnetization onto a sphere. 
It is demonstrated that the modified energy of the intermediate solution for both SAV schemes always decreases over time, which partially preserves energy dissipation.
Moreover, these two schemes exhibit high efficiency since they only require solving two linear systems with constant coefficient.
This combination of finite difference discretization and the Discrete Cosine Transform (DCT) leads to a highly efficient solution process, as the resulting linear systems can be solved using fast transform algorithms.
Numerical experiments demonstrate that both proposed SAV schemes possess significant advantages in terms of efficiency and fewer limitations on temporal step size for stability compared to common numerical methods.

The rest of the paper is organized as follows. 
In Section \ref{Sec:GF}, the micromagnetic flow model is introduced. 
Next, two SAV schemes are proposed for the Gibbs energy minimization problem in Section \ref{Sec:SAV}. 
And then, a series of numerical experiments are presented in Section \ref{Sec:NE} to  validate the effectiveness, energy dissipation, convergence, and efficiency of the proposed SAV schemes.
Finally, Section \ref{Sec:Conclusion} provides a conclusion.

%%%%%%%%%%%%%%%%%%%%%% Chapter Model %%%%%%%%%%%%%%%%%%%%%%%%%%%%%%%%%%%

\section{Micromagnetic flow model}\label{Sec:GF}

In this section, we present the Gibbs free energy in micromagnetics, then its corresponding gradient flow model is presented.

In micromagnetics, the state of magnetization is described by a vector-valued function $\h{M}(\h{x}, t) = (M_1(\h{x}, t), M_2(\h{x}, t), M_3(\h{x}, t))$, which maps a position $\h{x} \in \Omega$, where $\Omega \subset \mathbb{R}^3$ represents the region of the magnetic structure.
The length of the magnetization $\h{M}$ is consistently equal to the saturation magnetization $M_s$, namely $| \h{M}(\h{x}, t)  |= M_s$.
For simplicity, the normalized magnetization vector field is defined as $\h{m}(\h{x}, t) = \h{M}(\h{x}, t)/M_s$, which satisfies $|\h{m}(\h{x}, t)| \equiv 1$.

For a uniaxial magnetic material with an easy axis of $(1, 0, 0)$, the normalized free energy $g(\h{m})$, which consists of the exchange, anisotropy, magnetostatic, and Zeeman energies, is expressed as follows:
\begin{equation}\label{Eqn:FreeEnergy}
g(\h{m})
= \int_{\Omega}
\left( \frac{C_{ex}}{\mu_0 M_{\mathrm{s}}^2}|\nabla \h{m}|^2
+\frac{K_u}{\mu_0 M_{\mathrm{s}}^2}\left(m_2^2 + m_3^2 \right)
- \h{h}_e \cdot \h{m}
-\frac{1}{2} \h{h}_{{m}} \cdot \h{m}
\right) \mathrm{d} \h{x},
\end{equation}
where $C_{ex}$ is the exchange constant, 
$\mu_0$ is the magnetic permeability of vacuum,  
$K_u$ is the uniaxial anisotropy constant, 
$m_i$ denotes the $i$th component of $\h{m}$, 
$\h{h}_e$ is the external magnetic field, 
and $\h{h}_{{m}}$ is the magnetostatic field.

The effective field $\h{h}_{\mathrm{eff}}$, corresponding to the variational derivative of $g(\h{m})$ and comprising exchange, anisotropy, magnetostatic, and external magnetic fields, is represented as follows:
\begin{equation*}
\h{h}_{\mathrm{eff}} 
= - \dfrac{\delta g\left(\h{m}\right)}{\delta \h{m}}
= \frac{2 C_{ex}}{\mu_0 M_{\mathrm{s}}^2} \Delta \h{m} 
- \frac{2 K_u}{\mu_0 M_{\mathrm{s}}^2}  \left( m_2 \h{e}_2 + m_3 \h{e}_3 \right)
+ \h{h}_{{m}}
+ \h{h}_e, 
\end{equation*}
where $\h{e}_2 = (0, 1, 0)$ and $\h{e}_3 = (0, 0, 1)$.

Specially, the magnetostatic field $\h{h}_{{m}} = - \nabla U$ can be solved by (See \cite{YangLHuGH19:1048})
\begin{eqnarray}\label{Eqn:U}
\left\{\begin{aligned}
& \Delta U = \nabla \cdot \h{m}, && \text{in}\ \Omega, \\
& \Delta U = 0, && \text{outside}\ \Omega,\\
& [U]|_{\partial \Omega} = 0, 
\quad 
\left.\left[ \dfrac{\partial U} {\partial \h{n}}\right]\right|_{\partial \Omega} = - \h{m} \cdot \h{n}, 
\end{aligned}\right. 
\end{eqnarray}
where $\h{n}$ denotes the unit outward normal vector on $\partial \Omega$ and the jump $[U]_{\partial \Omega}$ is defined as
\begin{equation*}
[U]|_{\partial \Omega}(\h{x}) = \lim\limits_{\h{y} \rightarrow \h{x}, \h{y} \in \bar{\Omega}^c} U(\h{y}) - \lim\limits_{\h{y} \rightarrow \h{x},  \h{y} \in \Omega} U(\h{y}). 
\end{equation*}

The above Eq.~\eqref{Eqn:U} can be solved analytically, yielding
\begin{equation*}
U(\h{x}) =  \dfrac{1}{4 \pi} \int_\Omega \nabla_{\h{y}}\left( \dfrac{1}{|\h{x} - \h{y}|} \right) \cdot \h{m}(\h{y}) \mathrm{d} \h{y}. 
\end{equation*}
Consequently, the analytical solution of the magnetostatic field $\h{h}_{{m}} = - \nabla U$ is written as
\begin{equation}\label{Eqn:hm}
\h{h}_{{m}} = - \dfrac{1}{4 \pi}  \nabla_{\h{x}} \int_\Omega \nabla_{\h{y}}\left( \dfrac{1}{|\h{x} - \h{y}|} \right) \cdot \h{m}(\h{y}) \mathrm{d} \h{y}. 
\end{equation}

\begin{remark}\label{Rem:MagnetostaticEenergy}
By Eq.~\eqref{Eqn:U}, the magnetostatic energy ${E}_m$ can be equivalently expressed as:
\begin{equation*}
{E}_m = -\frac{1}{2} \int_\Omega \h{h}_{{m}} \cdot \h{m} \mathrm{d} \h{x}
= \frac{1}{2}\int_{\mathbb{R}^3} |\nabla U|^2 \mathrm{d} \h{x} \geq 0. 
\end{equation*}
\end{remark}

To obtain the magnetic ground state, we consider the Gibbs energy minimization problem:
\begin{equation}\label{Eqn:MEMP}
\min\limits_{\h{m}} g(\h{m}) \quad \text{with constraint}\ |\h{m}| = 1, 
\end{equation}
which constitutes a non-convex and nonlinear optimization problem due to the pointwise constraint.

The flow equation of the Gibbs energy minimization problem \eqref{Eqn:MEMP} can be given directly as 
\begin{equation}\label{Eqn:GF}
\dfrac{\partial \h{m}}{\partial t} = \h{h}_{\mathrm{eff}},  
\end{equation}
with the pointwise constraint $|\h{m}| = 1$.
This equation is subject to the homogeneous Newman boundary condition  ${\partial \h{m}}/{ \partial \h{n}} = \h{0}$.

Governed by an energy dissipation law, the flow equation \eqref{Eqn:GF} ensures that the free energy $g(\h{m})$ decreases over time when the external magnetic field is held constant, thus describing an energy-dissipative process.
Forming the discrete inner product with $\h{m}_t$ on both sides of Eq.~\eqref{Eqn:GF}, we can directly obtain
\begin{equation}\label{Eqn:LS}
\dfrac{\mathrm{d} g(\h{m})}{\mathrm{d} t} = - \int_\Omega |\h{m}_t|^2 \mathrm{d}\h{x}. 
\end{equation}
Eq.~\eqref{Eqn:LS} guarantees that the system tends toward stable equilibrium, which is a crucial characteristic for energy minimization problems.

In the remainder of this paper, we consistently assume that the external magnetic field $\h{h}_e$
is zero.

\section{Numerical schemes}\label{Sec:SAV}

In this section, some common numerical schemes are first presented to motivate our methods.
Next, two SAV schemes are proposed, along with the analyses of energy dissipation and computational cost.
It is noted that we focus on temporal discretization and do not specifically address spatial discretization in this section.

\subsection{Numerical schemes for flow equation}\label{Sec:BCFE}

For the temporal discretization, a first-order implicit projection method was proposed for the Landau-Lifshitz (LL) equation when concentrating on the exchange field in \cite{EWNWangXP01:1647}. 
This method possesses the property of unconditional stability.
First, we present a specific numerical scheme as follows: 
\begin{equation}\label{Eqn:PM}
\left\{\begin{aligned}
&\dfrac{\h{m}^* - \h{m}^n}{\Delta t} = \Delta \h{m}^*,  
\\
&\h{m}^{n+1} = \dfrac{\h{m}^*}{|\h{m}^*|}. 
\end{aligned}\right. 
\end{equation}
Here, $\h{m}^n$ denotes the magnetization at time $t_n = n \Delta t$, with $\Delta t$ denoting the temporal step size.
Algorithmically, the projection method \eqref{Eqn:PM} first solves the linear flow equation \eqref{Eqn:GF} and then projects the solution onto the unit sphere.

While the efficiency of solving linear systems makes the projection method an attractive candidate for full Gibbs energy minimization, a direct application proves computationally inefficient. To establish a baseline and highlight the need for our approach, we first evaluate classical time-stepping schemes (e.g., forward and backward Euler). These methods are shown to be inadequate, as they often violate the energy dissipation law and incur high computational costs.

To contextualize our proposed method, we first examine the limitations of two standard Euler-type projection schemes for solving the micromagnetic flow equation.

The Forward Euler Projection (FEP) scheme, defined by
\begin{equation}\label{Eqn:FE}
\left\{\begin{aligned}
&\dfrac{\h{m}^{*} - \h{m}^{n}}{\Delta t} 
=   \h{h}_{\mathrm{eff}}(\h{m}^{n}), 
\\
 &\h{m}^{n+1} = \dfrac{\h{m}^{*}}{|\h{m}^{*}|}, 
\end{aligned}\right.
\end{equation}
offers computational simplicity. However, its explicit nature fundamentally violates the energy dissipation structure, often requiring restrictive time-step constraints to maintain numerical stability, as detailed in Section \ref{Sec:Efficiency}.

In contrast, the Backward Euler Projection (BEP) scheme,
\begin{equation}\label{Eqn:BE}
\left\{\begin{aligned}
&\dfrac{\h{m}^{*} - \h{m}^{n}}{\Delta t} 
=   \h{h}_{\mathrm{eff}}(\h{m}^{*}), 
\\
 &\h{m}^{n+1} = \dfrac{\h{m}^{*}}{|\h{m}^{*}|}. 
\end{aligned}\right.
\end{equation}
 
The BEP scheme is a simple extension of projection method \eqref{Eqn:PM}.
The intermediate value $\h{m}^{*}$ satisfies an approximate energy dissipation law.
The energy dissipation law is derived by first taking the inner product of the first equation in Eq.~\eqref{Eqn:BE} with $({\h{m}^{*} - \h{m}^{n}})/{\Delta t}$. Subsequently, application of Green's formula, $2a (a-b) = a^2 - b^2 + (a-b)(a-b)$, yields the following result:
\begin{eqnarray*}
\dfrac{g(\h{m}^{*}) - g(\h{m}^{n})}{\Delta t} &=&
-  \left\| \dfrac{\h{m}^{*} - \h{m}^{n}}{\Delta t}  \right\|^2 
- \dfrac{1}{\Delta t}\dfrac{C_{ex}}{\mu_0 M_s^2} \| \nabla (\h{m}^{*} - \h{m}^{n})\|^2 
\\
&& 
 + \dfrac{1}{2 \Delta t}  \Big(\h{h}_m(\h{m}^{*}) - \h{h}_m(\h{m}^{n}), \h{m}^{*} - \h{m}^{n}\Big)
\\
&&
- \dfrac{1}{\Delta t}\dfrac{K_u}{\mu_0 M_s^2} \Big( \| {m}^{*}_2 - {m}^{n}_2\|^2 +  \| {m}^{*}_3 - {m}^{n}_3\|^2\Big),  
\end{eqnarray*}
where $(\cdot, \cdot)$ and $\| \cdot \|$ denote the $L^2$ inner product and $L^2$ norm, respectively.
By Remark \ref{Rem:MagnetostaticEenergy}, it can be seen that the BEP scheme possesses a discrete energy dissipation law.

Unfortunately, the implicit scheme of the magnetostatic field $\h{h}_m(\h{m}^{*})$ leads to an unacceptably high computational expense, as the magnetostatic field needs to be integrated over the entire computational domain $\Omega$ as shown in Eq.~\eqref{Eqn:hm}.
In the implementation, the non-local property of the magnetostatic field results in a dense matrix after spacial discretization. 
Hence, an iterative method for the BEP scheme is essential, which will result in a substantial increase in computational time.
Therefore, although the BEP scheme can preserve energy dissipation unconditionally, it suffers from significant computational time.

\subsection{SAV discrete schemes}

By introducing scalar auxiliary variables, the SAV method not only ensures energy dissipation, but also can be solved efficiently.
To address the challenges of both FEP and BEP schemes, two SAV schemes are proposed to solve the flow equation \eqref{Eqn:GF}.

To mitigate the high computational cost associated with the magnetostatic field, we introduce a scalar auxiliary variable $r$ for the magnetostatic energy, which is defined by 
\begin{equation}\label{Eqn:r}
r = \sqrt{- \dfrac{1}{2}\int_\Omega \h{h}_m \cdot \h{m} \mathrm{d} x}. 
\end{equation}
The scalar auxiliary variable $r$ is always well-defined due to the non-negativity of the magnetostatic energy (See Remark \ref{Rem:MagnetostaticEenergy}).

For simplicity, we denote the following physical constants as
$$
C_{e} :=  \dfrac{2 C_{ex}}{\mu_0 M_s^2}   
\quad 
\text{and}
\quad
C_{an} : = \dfrac{2 K_u}{\mu_0 M_s^2}. 
$$

By introducing the scalar auxiliary variable $r$, Eq.~\eqref{Eqn:GF} can be reformulated and the SAV method is proposed as follows:
\begin{subequations}\label{Eqn:SAV1}
\begin{align}\label{Eqn:SAV1a}
& \h{m}_t =
 C_{e}\Delta \h{m}
- C_{an} (m_2 \h{e}_2 + m_3 \h{e}_3)
+ \dfrac{r}{\sqrt{- \frac{1}{2}( \h{h}_m , \h{m} )}}   \h{h}_m, 
\\
& \label{Eqn:SAV1b}
  \dfrac{\mathrm{d} r}{\mathrm{d} t}  = \left( \frac{\delta r}{\delta \h{m}} , \h{m}_t \right) = -  \dfrac{1}{2 r}  (  \h{h}_m, \h{m}_t) . 
\end{align}
\end{subequations}

The SAV method \eqref{Eqn:SAV1} automatically preserves energy dissipation.
Indeed, by taking inner product of both sides of Eq.~\eqref{Eqn:SAV1a} with $\h{m}_t$ and multiplying both sides of Eq.~\eqref{Eqn:SAV1b} with $2r$, we can obtain that 
\begin{eqnarray}\nonumber
 \left\| \h{m}_t\right\|^2 &=&
C_{e} \left(\Delta \h{m}, \h{m}_t\right)
- C_{an} \Big[ (m_2, m_{t, 2}) +  (m_3, m_{t, 3})\Big]
\\ \label{Eqn:SAV1Lin}
&&+\frac{r}{\sqrt{-\frac{1}{2} ( \h{h}_m, \h{m} )}}  ( \h{h}_m, \h{m}_t), 
\\  \label{Eqn:SAV2Lin}
 2 r r_t &=& - ( \h{h}_m, \h{m}_t), 
\end{eqnarray}
where $m_{t, i}$ denotes the $i$th component of $\h{m}_t$. 
%%%
%%%
By substituting Eq.~\eqref{Eqn:SAV2Lin} into Eq.~\eqref{Eqn:SAV1Lin} and applying  Green's formula with boundary condition $\partial \h{m}/ \partial \h{n} = \h{0}$ and Eq.~\eqref{Eqn:r}, we obtain
$$
\dfrac{\mathrm{d} \tilde{g}(\h{m}, r)}{\mathrm{d} t} = -  \|\h{m}_t\|^2,
$$
where $\tilde{g}(\h{m})$ is the modified energy defined as
\begin{equation}\label{Eqn:ModifiedEnergy}
\tilde{g}(\h{m}, r) = 
\dfrac{C_{e}}{2}  \int_\Omega |\nabla \h{m}|^2\mathrm{d} \h{x} 
+ \dfrac{C_{an}}{2} \int_\Omega   \left( m_2^2 + m_3^2 \right) \mathrm{d} \h{x} 
+ r^2. 
\end{equation}
Therefore, the SAV method \eqref{Eqn:SAV1} guarantees the energy dissipation in the sense of the modified energy $\tilde{g}(\h{m}, r)$.

{\bf  SAV discrete scheme.} We present a first-order temporal discrete scheme SAV1 of the SAV method \eqref{Eqn:SAV1} as follows:
\begin{subequations}\label{Eqn:SAVTD}
\begin{align}\label{Eqn:SAVTDa}
&\frac{\h{m}^{*}- \h{m}^n}{\Delta t} =
C_{e} \Delta_h \h{m}^{*}
- C_{an} ({m}^{*}_2 \h{e}_{2} + {m}^{*}_3 \h{e}_{3})
+\frac{r^{*}}{\sqrt{-\frac{1}{2} (\h{h}_m^n, \h{m}^n)_h}}   \h{h}_m^n,
\\ \label{Eqn:SAVTDb}
& \dfrac{r^{*}-r^n}{\Delta t} = -\frac{1}{2 \sqrt{-\frac{1}{2} (\h{h}_m^n, \h{m}^n)_h}} \left(\h{h}_m^n,  \dfrac{\h{m}^{*}- \h{m}^n}{\Delta t} \right)_h, 
\\ \label{Eqn:SAVTDc}
 & \h{m}^{n+1}= \dfrac{\h{m}^{*}}{|\h{m}^{*}|}, 
\\ \label{Eqn:SAVTDd}
 &\dfrac{r^{n+1}-r^n}{\Delta t} = -\frac{1}{2 \sqrt{-\frac{1}{2} (\h{h}_m^n, \h{m}^n)_h}} \left(\h{h}_m^n,  \dfrac{\h{m}^{n+1}- \h{m}^n}{\Delta t} \right)_h, 
\end{align}
\end{subequations}
where $\Delta_h $ denotes the discrete Laplace operator, $(\cdot, \cdot)_h$ denotes the discrete $L^2$ inner product, and $\h{h}_m^n$ denotes $\h{h}_m$ at time $t_n = n \Delta t$.  
The corresponding discrete $L^2$ norm is denoted as $\| \h{u}\|_h = \sqrt{(\h{u}, \h{u})_h}$. 

It can be seen that Eqs.~\eqref{Eqn:SAVTDa} and \eqref{Eqn:SAVTDb} are the discretizations of Eqs.~\eqref{Eqn:SAV1a} and \eqref{Eqn:SAV1b}, respectively. 
Noticing that Eqs.~\eqref{Eqn:SAVTDa} and \eqref{Eqn:SAVTDb} cannot guarantee $|\h{m}^*| = 1$, the simple projection is employed to obtain $\h{m}^{n+1}$ in Eq.~\eqref{Eqn:SAVTDc} to preserve the constant magnetization length. 
To avoid the errors caused by an incorrect magnetization length, $r^{n+1}$ is updated by using Eq.~\eqref{Eqn:SAVTDd} with $\h{m}^{n+1}$ rather than Eq.~\eqref{Eqn:SAVTDb}.

From the numerical experiments in Section \ref{Sec:NE}, we can observe a small discrepancy between the scalar auxiliary variable $r$ and the magnetostatic energy when employing a large temporal step size.
To enhance the correspondence between $r$ and the magnetostatic energy, we modify the update of $r$ based on its definition and propose the following modified SAV discrete scheme (A similar idea can be found in \cite{ZhuangQQShenJ19:72}).

{\bf Modified SAV discrete scheme.} The first-order modified SAV scheme SAV2 is presented as follows:  
\begin{subequations}\label{Eqn:RTUT}
\begin{align}\label{Eqn:RTUTa}
&\frac{\h{m}^{*}- \h{m}^n}{\Delta t} =
C_{e} \Delta_h \h{m}^{*}
- C_{an} ({m}^{*}_2 \h{e}_{2} + {m}^{*}_3 \h{e}_{3})
+\frac{r^{*}}{\sqrt{-\frac{1}{2} (\h{h}_m^n, \h{m}^n)_h}}   \h{h}_m^n,
\\ \label{Eqn:RTUTb}
& \dfrac{r^{*}-r^n}{\Delta t} = -\frac{1}{2 \sqrt{-\frac{1}{2} (\h{h}_m^n, \h{m}^n)_h}} \left(\h{h}_m^n,  \dfrac{\h{m}^{*}- \h{m}^n}{\Delta t} \right)_h, 
\\ \label{Eqn:RTUTc}
 & \h{m}^{n+1}= \dfrac{\h{m}^{*}}{|\h{m}^{*}|}, 
\\ \label{Eqn:RTUTd}
 &r^{n+1} = \sqrt{- \frac{1}{2}(\h{h}_m^{n+1}, \h{m}^{n+1})_h }.  
\end{align}
\end{subequations}

Comparing with SAV1 scheme, SAV2 scheme updates $r$ by the definition of $r$ instead of solving dynamic equation Eq.~\eqref{Eqn:SAVTDd}, in which the correspondence between $r$ and the magnetostatic energy can be effectively improved.
Indeed, $r^2$ is exactly equal to the magnetostatic energy in this case, which implies that the original energy and the modified energy are the same.

\subsection{Energy dissipation of discrete SAV schemes}

In the following, we will analyze the energy dissipation of both SAV1 scheme \eqref{Eqn:SAVTD} and SAV2 scheme \eqref{Eqn:RTUT}.

By taking  discrete $L^2$ inner product with $\frac{\h{m}^{*} - \h{m}^n}{\Delta t}$ on both sides of Eq.~\eqref{Eqn:SAVTDa} and multiplying both sides of Eq.~\eqref{Eqn:SAVTDb} by $2r^{*}$, we obtain
\begin{eqnarray} \nonumber
 \left\| \frac{\h{m}^{*}- \h{m}^n}{\Delta t} \right\|^2_h 
  &=&
\dfrac{C_{e}}{\Delta t} (\Delta_h \h{m}^{*}, \h{m}^{*}- \h{m}^n)_h
\\ \nonumber
&&
-\dfrac{C_{an}}{\Delta t} \Big[ ({m}^{*}_2, {m}^{*}_2 - {m}^{n}_2)_h + ({m}^{*}_3, {m}^{*}_3 - {m}^{n}_3)_h\Big]
\\ \label{Eqn:SAVTDamn+1-mn}
&&
+\dfrac{1}{\Delta t}\dfrac{r^{*}}{\sqrt{-\frac{1}{2} (\h{h}_m^n, \h{m}^n)_h}}   (\h{h}_m^n, \h{m}^{*}- \h{m}^n)_h,
\\  \label{Eqn:SAVTDbr}
 2 (r^{*}-r^n)r^{*} &=& -\frac{r^{*}}{ \sqrt{-\frac{1}{2} (\h{h}_m^n, \h{m}^n)_h}} \left( \h{h}_m^n,  \h{m}^{*}- \h{m}^n \right)_h. 
\end{eqnarray}
%%%%
%%%%
By substituting Eq.~\eqref{Eqn:SAVTDbr} into Eq.~\eqref{Eqn:SAVTDamn+1-mn} and applying the identity $2a (a-b) = a^2 - b^2 + (a-b)(a-b)$, we can derive that
\begin{eqnarray} \nonumber
\frac{\tilde{g}_h(\h{m}^{*}, r^*)-\tilde{g}_h (\h{m}^{n}, r^n)}{\Delta t} 
&=&
- \left\|\frac{\h{m}^{*}-\h{m}^n}{\Delta t}\right\|_h^2
-\frac{C_{an}}{2 \Delta t} \Big( \left\|{m}^{*}_2- {m}^n_2\right\|_h^2 + \left\|{m}^{*}_3- {m}^n_3\right\|_h^2\Big)
\\ \nonumber
&& 
- \frac{1}{\Delta t} \left(r^{*}-r^n\right)^2
+ \dfrac{C_{e}}{2\Delta t}  \Big(\Delta_h (\h{m}^{*} - \h{m}^{n}), \h{m}^{*} - \h{m}^{n} \Big)_h
\\ \label{Eqn:DiscreteEnergy}
&& + \dfrac{C_{e}}{2\Delta t}  \Big[ (\Delta_h \h{m}^n, \h{m}^{*})_h -  (\Delta_h \h{m}^{*}, \h{m}^n)_h\Big], 
\end{eqnarray}
where the discrete modified energy $\tilde{g}_h$ is defined as
$$
\tilde{g}_h (\h{m}^{n}, r^n) =
-  \dfrac{C_{e}}{2} (\Delta_h \h{m}^{n}, \h{m}^{n})_h 
+ \dfrac{C_{an}}{2} \Big(  \left\|{m}^n_2 \right\|_h^2 + \left\|{m}^n_3 \right\|_h^2 \Big)
+\left(r^n\right)^2. 
$$

By Eq.~\eqref{Eqn:DiscreteEnergy}, the following theorem is established.
\begin{theorem}
Assume $\h{m}^{*}$ is obtained by SAV1 scheme \eqref{Eqn:SAVTD} or SAV2 scheme \eqref{Eqn:RTUT}.
If the following two conditions
\begin{eqnarray}\label{Eqn:Con1}
& (\Delta_h (\h{m}^{*} - \h{m}^{n}), \h{m}^{*} - \h{m}^{n})_h  \leq 0, 
\\ \label{Eqn:Con2}
& (\Delta_h \h{m}^n, \h{m}^{*})_h -  (\Delta_h \h{m}^{*}, \h{m}^n)_h  = 0, 
\end{eqnarray}
are satisfied, then the energy dissipation inequality
\begin{eqnarray*} \nonumber
\tilde{g}_h(\h{m}^{*}, r^*)
\leq \tilde{g}_h (\h{m}^{n}, r^n) 
\end{eqnarray*}
holds. 
\end{theorem}

\begin{remark}
For finite element spatial discretization, the conditions \eqref{Eqn:Con1} and \eqref{Eqn:Con2} can be achieved by directly employing the Green's formula with the boundary condition $\partial \h{m}/\partial \h{n} = 0$.
In the context of finite difference spatial discretization, these two conditions can also be satisfied by employing the second-order central difference for the discrete Laplace operator (More details are refered to Lemma 2.3 in \cite{ChenJRWangC21:55}). 
\end{remark}

\begin{remark}\label{Rem:Energy}
Because of the simple projection of $\h{m}^*$ and the update of $r^{n+1}$ in both SAV1 and SAV2 schemes, the strong result $\tilde{g}_h(\h{m}^{n+1}, r^{n+1}) \leq \tilde{g}_h (\h{m}^{n}, r^{n})$ cannot be proved theoretically. 
Only an intermediate result $\tilde{g}_h(\h{m}^{*}, r^*) \leq \tilde{g}_h (\h{m}^{n}, r^n)$ is obtained. 
However, numerical experiments demonstrate that the modified energy $\tilde{g}_h(\h{m}^{n+1}, r^{n+1})$ remains dissipative during the magnetization evolution towards the equilibrium state. 
For more details, we refer the readers to Section \ref{Sec:effectiveness}. 
\end{remark}

\subsection{Analysis of computational cost}

In addition to preserving energy dissipation, another remarkable advantage of SAV method is its exceptional computational efficiency. 
In this part, we will analyze the computational cost associated with the two proposed SAV schemes.

By Eq.~\eqref{Eqn:SAVTDb}, we have
\begin{eqnarray}\label{Eqn:rn+1}
 r^{*} = r^n-\frac{1}{2 \sqrt{-\frac{1}{2} (\h{h}_m^n, \h{m}^n)_h}} (\h{h}_m^n,  \h{m}^{*}- \h{m}^n)_h. 
\end{eqnarray}
Substituting the above equation into Eq.~\eqref{Eqn:SAVTDa}, we obtain
\begin{eqnarray}\label{Eqn:Rewritten}
 A \h{m}^{*} - \Delta t \frac{(\h{h}_m^n,  \h{m}^{*})_h}{ (\h{h}_m^n, \h{m}^n)_h}  \h{h}_m^n
=   \h{F},
\end{eqnarray}
where the operator $A$ and notation $\h{F}$ are defined as
\begin{eqnarray*}\label{Eqn:A}
 &A \h{m}^{*} = \h{m}^{*} 
+  C_{an} \Delta t ({m}^{*}_2  \h{e}_{2} + {m}^{*}_3  \h{e}_{3} )
- C_{e}  \Delta t   \Delta_h \h{m}^{*} , 
\\
 & \h{F}= \h{m}^n +  \Delta t 
\left( \dfrac{r^n}{\sqrt{-\frac{1}{2} (\h{h}_m^n, \h{m}^n)_h}} -1 \right)  \h{h}_m^n, 
\end{eqnarray*}
respectively. 
%%%%%%%
%%%%%%%
Left-multiplication of both sides of Eq.~\eqref{Eqn:Rewritten} by $A^{-1}$ yields
\begin{eqnarray}\label{Eqn:Mn+1A-1}
\h{m}^{*}
- \Delta t \frac{(\h{h}_m^n,  \h{m}^{*})_h}{  (\h{h}_m^n, \h{m}^n)_h}  A^{-1} \h{h}_m^n
=  A^{-1} \h{F}. 
\end{eqnarray}
Taking the discrete inner product of both sides of the preceding equation against $\h{h}_m^n$ leads to
\begin{eqnarray}\nonumber
(\h{m}^{*},  \h{h}_m^n)_h
- \Delta t \frac{(\h{h}_m^n,  \h{m}^{*})_h}{  (\h{h}_m^n, \h{m}^n)_h}  (A^{-1} \h{h}_m^n, \h{h}_m^n)_h
=  (A^{-1} \h{F}, \h{h}_m^n)_h,
\end{eqnarray}
which could be derived to
\begin{eqnarray}\label{Eqn:Mn+1Hmn}
(\h{m}^{*},  \h{h}_m^n)_h =  (A^{-1} \h{F}, \h{h}_m^n)_h \left/ \left(1 - \Delta t \frac{(A^{-1} \h{h}_m^n, \h{h}_m^n)_h}{ (\h{h}_m^n, \h{m}^n)_h}   \right)\right..
\end{eqnarray}

In the implement process of SAV1 and SAV2 schemes, we first solve Eq.~\eqref{Eqn:Mn+1Hmn} to obtain $(\h{m}^{*},  \h{h}_m^n)_h$. 
This step involves solving two linear systems with constant coefficient, namely solving $A\h{x} = \h{F}$ and $A\h{y} =  \h{h}_m^{n}$, and obtaining $A^{-1}\h{F}$ and $A^{-1}\h{h}_m^{n}$.
And then, $\h{m}^{*}$ can be directly obtained by Eq.~\eqref{Eqn:Mn+1A-1}.
At last, $\h{m}^{n+1}$ and $r^{n+1}$ are accordingly updated by Eqs.~\eqref{Eqn:SAVTDc} and \eqref{Eqn:SAVTDd} (or \eqref{Eqn:RTUTd}), respectively. 
%%%%
Consequently, the primary computational cost arises from solving two linear systems with constant coefficient.

\begin{remark}
As $r^{n+1}$ is updated by Eq.~\eqref{Eqn:SAVTDd} or \eqref{Eqn:RTUTd}, the solving for $r^*$ in Eq.~\eqref{Eqn:rn+1} is not necessary.
\end{remark}

\begin{remark}
The solutions of the two linear systems $A\h{x} = \h{F}$ and $A\h{y} =  \h{h}_m^{n}$ can be efficiently obtained by using the Discrete Cosine Transform when finite difference spatial discretization is employed.
\end{remark}

\section{Numerical experiments}\label{Sec:NE}

This section is devoted to numerical validation of the SAV1 \eqref{Eqn:SAVTD} and SAV2 \eqref{Eqn:RTUT} schemes, assessing their effectiveness, energy dissipation, convergence, and computational efficiency.
It it noted that both finite difference and finite element spatial discretizations are applicable to the proposed SAV schemes.
In our implementation, the finite different spacial discretization is employed with a discrete Laplace operator of second-order center difference. 
Additionally, the computation of the magnetostatic field $\h{h}_{m}$ is achieved by using 3D Fast Fourier Transform. 
Alternative methods include the fast multipole method \cite{BlueJLScheinfeinMR91:4778, CuiZYangL24:}, the tensor grid method \cite{ExlLAuzingerW12:2840}, and others.
All numerical experiments are conducted using a machine with Intel(R) Core(TM) i7-10750H CPU 2.60 GHz and 16-GB memory.

\subsection{Effectiveness and convergence}\label{Sec:effectiveness}

\begin{figure}[!t]
\begin{center}
\subfigbottomskip = 1pt % set the distant between the figure in a row
\subfigcapskip = -2pt % set the distant between the subfigure and subtitile
\subfigure[Initial magnetization]{\label{Fig:SAVMRTUF-Initial-M}
\includegraphics[width=69mm]{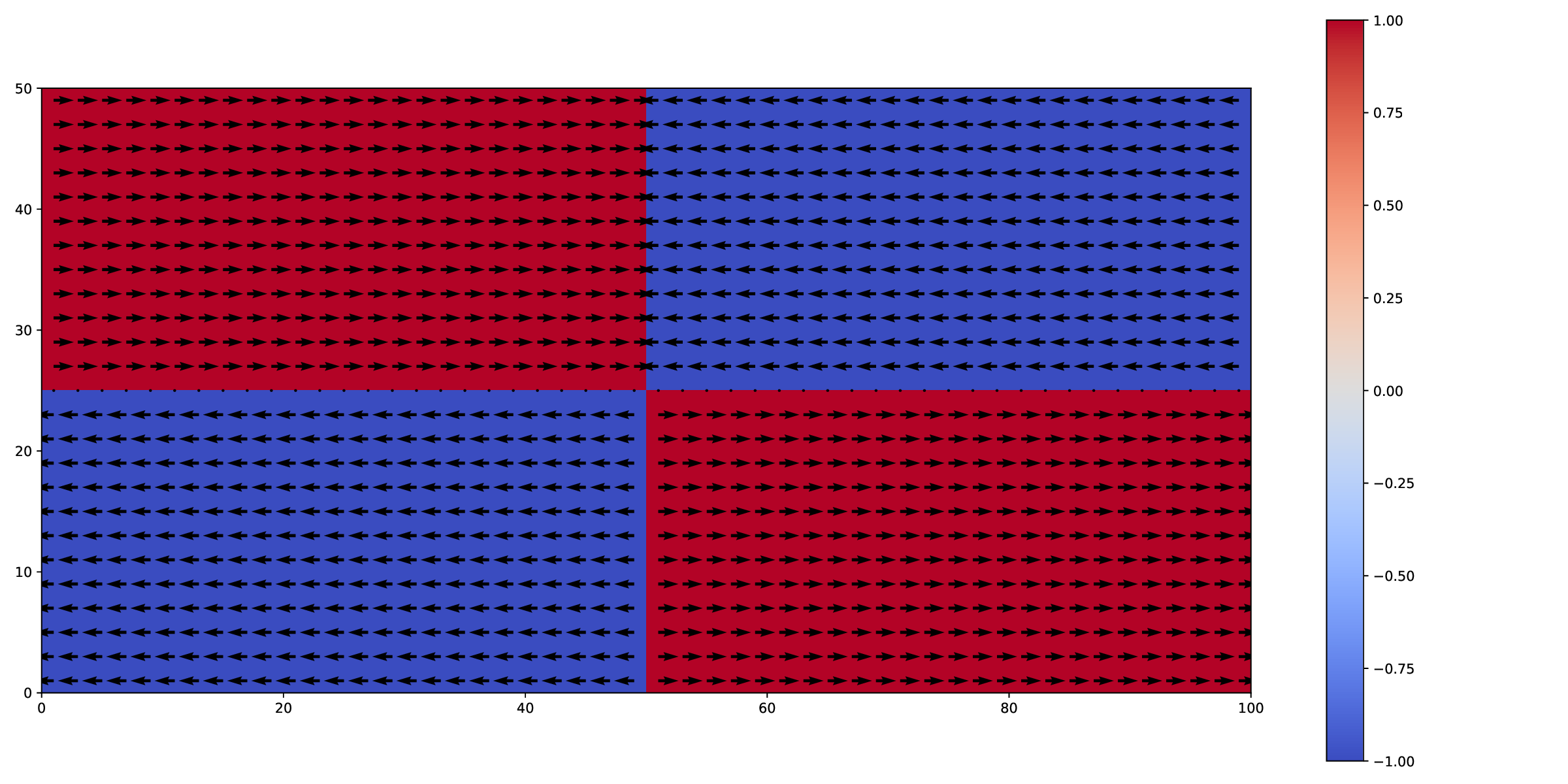}  %  fig 表示路径
}
\subfigure[Diamond state]{\label{Fig:SAVMRTUF-Diamond-M}
\includegraphics[width=69mm]{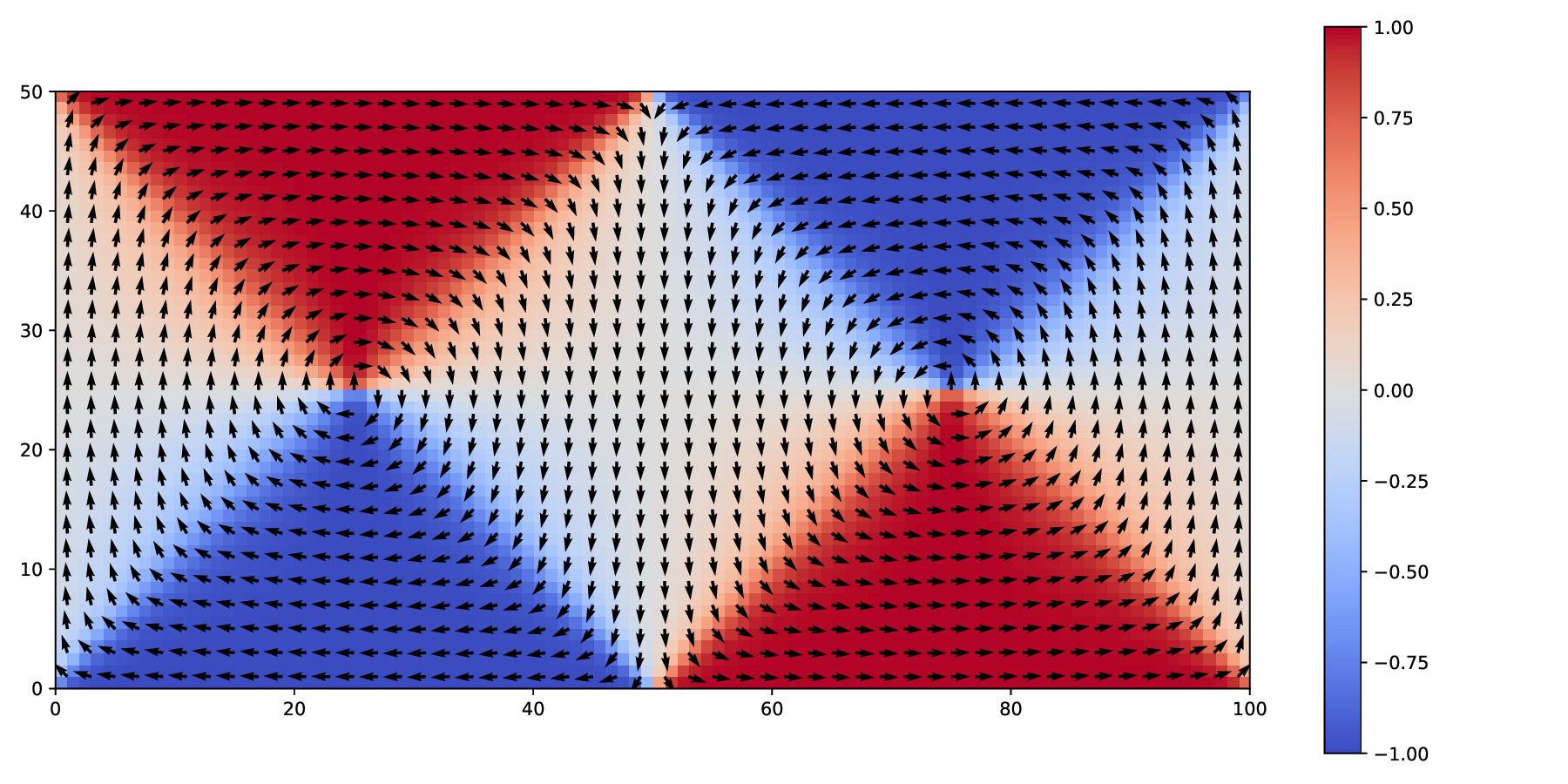}  %  fig 表示路径
}
\caption{(a) Initial magnetization \eqref{Eqn:IMDiamond}; (b) the corresponding steady-state named as diamond state. The color mapping indicates the angle between the magnetization vector and the $x$-axis.}
%The background color represents the angle between the magnetization and $x$-axis. }
\label{Fig:IMDiamond}
\end{center}
\end{figure}

We test the effectiveness of SAV1 and SAV2 schemes by the following example.
\begin{example}\label{Exa:1}
Let $\Omega$ be a ferromagnetic thin film of size $2 \mu\mathrm{m} \times 1 \mu\mathrm{m} \times 0.02 \mu\mathrm{m}$. 
The physical parameters are identical to those in Standard Problem No.1 \cite{muMag} as: $\alpha = 0.1$, $\gamma = 2.211 \times 10^5$ m/(As), $C_{ex} = 1.3 \times 10^{-11}$ J/m, 
$\mu_0 = 4\pi \times 10^{-7}$ N/$A^2$,
$M_s = 8.0 \times 10^5$ A/m, 
$K_u = 5.0 \times 10^2$ J/$m^3$. 
%$\h{e}_{an}  = (1, 0, 0)$. 
%
Initial magnetization (unit of length is $\mu \mathrm{m}$) are given as:
\begin{equation}\label{Eqn:IMDiamond}
\h{m}^0(x, y, z) = \left\{\begin{aligned}
& (-1,0, 0), && \mbox{if $0  \leq x \leq 1 , 0  \leq y \leq 0.5 $}, \\
& (1,0, 0), && \mbox{if $1 \leq x \leq 2 , 0  \leq y \leq 0.5 $}, \\
& (1,0, 0), && \mbox{if $0  \leq x \leq 1 , 0.5 \leq y \leq 1 $}, \\
& (-1,0, 0), && \mbox{if $1 \leq x \leq 2 , 0.5  \leq y \leq 1 $}, 
\end{aligned}\right.
\end{equation}
which can induce a diamond state as the steady state of magnetization \cite{RaveWHubertA00:3886}, as shown in Fig.~\ref{Fig:IMDiamond}.
\end{example}
%%%

Example \ref{Exa:1} is implemented using a uniform spacial mesh of $(100, 50, 1)$.
To obtain a comparable temporal step size to that of the micromagnetic dynamic equations, we adopt a time rescaling $t \rightarrow \eta t$, where $\eta = \alpha/(\gamma M_s)$ is a positive phenomenological damping parameter. 
This time rescaling is applied to the flow equation \eqref{Eqn:GF}, resulting in the following equation
\begin{equation}\label{Eqn:etaGF}
\eta \h{m}_{t} = \h{h}_{\mathrm{eff}}. 
\end{equation}

\begin{figure}[!t]
\begin{center}
\subfigbottomskip = 1pt % set the distant between the figure in a row
\subfigcapskip = -2pt % set the distant between the subfigure and subtitile
\subfigure[Original and modified total energies]{\label{Fig:SAVMRTUF-Diamond-Energy}
\includegraphics[width=69mm]{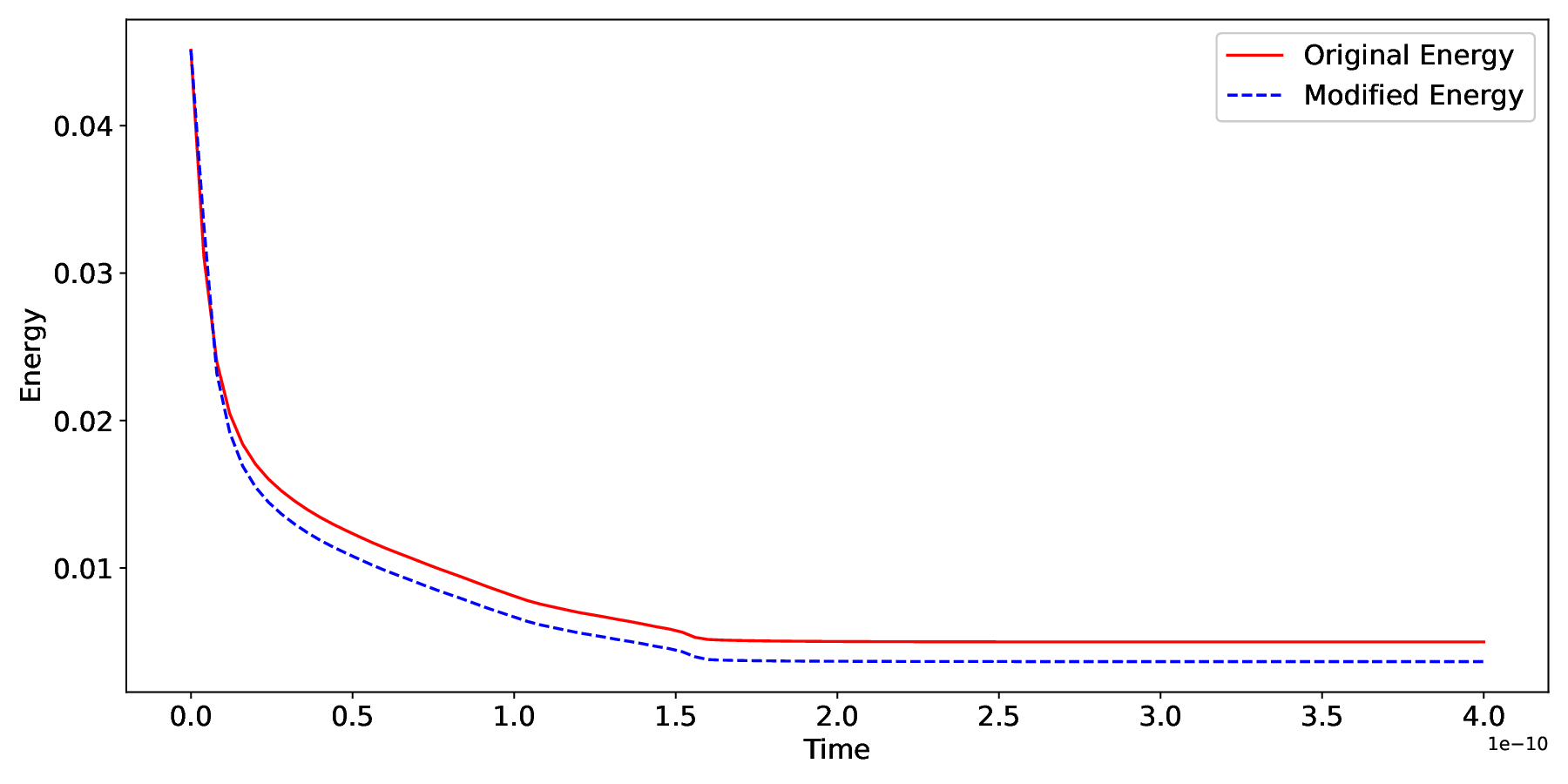}  %  fig 表示路径
}
\subfigure[Original and modified magnetostatic energy]{\label{Fig:SAVMRTUF-Diamond-SFEnergy}
\includegraphics[width=69mm]{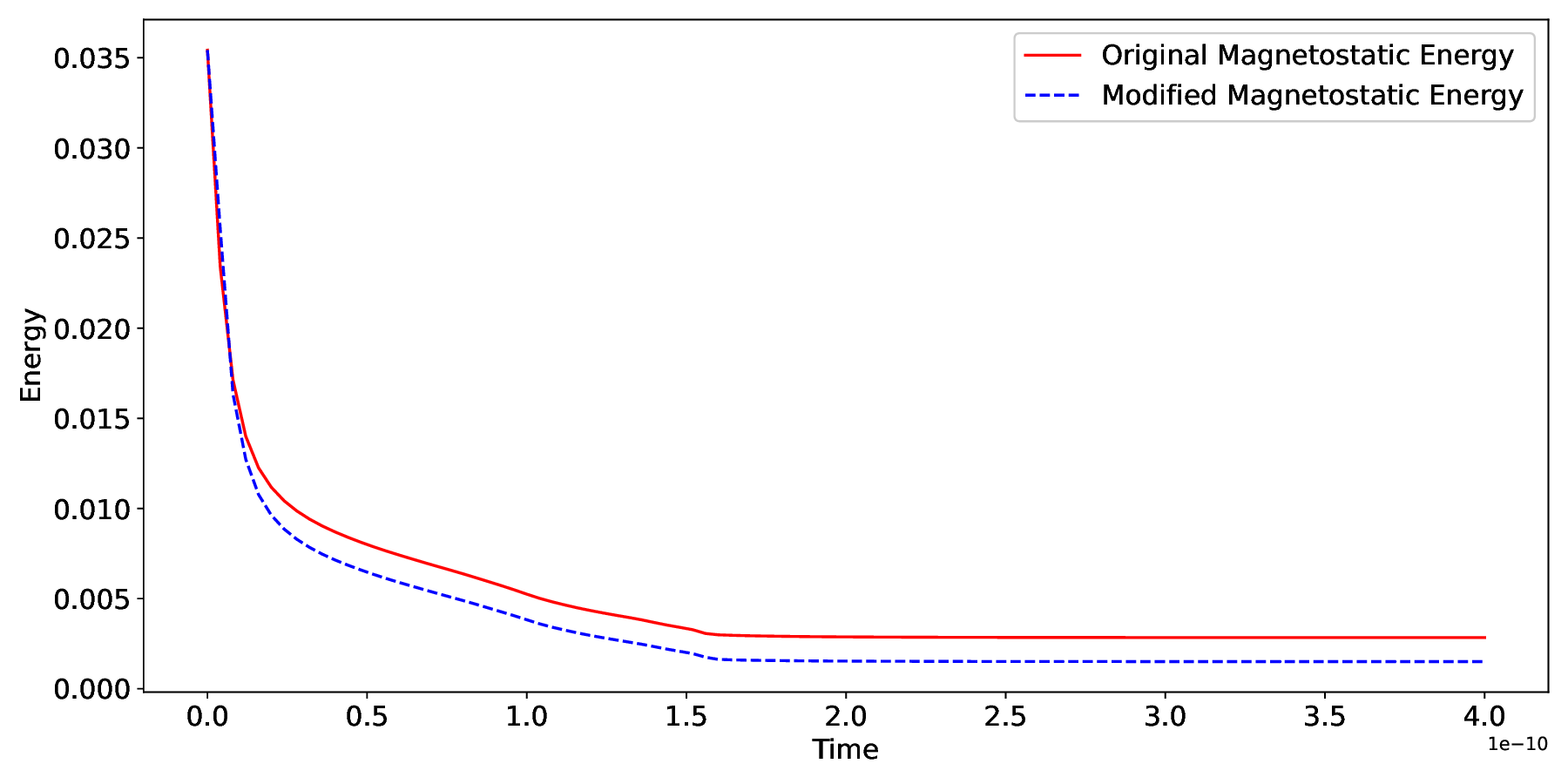}  %  fig 表示路径
}
\caption{Comparison of energies as functions of time by applying SAV1 scheme with initial magnetization \eqref{Eqn:IMDiamond} and $\Delta t = 10^{-12}$.}
\label{Fig:SAVMRTUF-Diamond}
\end{center}
\end{figure}

(i) Firstly, the effectiveness of the SAV1 scheme \eqref{Eqn:SAVTD} is validated with temporal step size $\Delta t = 10^{-12}$ and total simulation time $T = 4\times 10^{-10}$.
From Fig.~\ref{Fig:SAVMRTUF-Diamond-Energy}, it can be observed that both the original energy \eqref{Eqn:FreeEnergy} and modified energy \eqref{Eqn:ModifiedEnergy} smoothly and steadily decrease, eventually an equilibrium state is reached.
However, there is a small discrepancy between the original and modified energy.
The original magnetostatic energy and modified magnetostatic energy $r^2$, defined in Eqs.~\eqref{Eqn:FreeEnergy} and \eqref{Eqn:ModifiedEnergy}, respectively, are presented as functions of time in Fig.\ref{Fig:SAVMRTUF-Diamond-SFEnergy}.
From Fig.~\ref{Fig:SAVMRTUF-Diamond-SFEnergy}, it can be confirmed that the difference between the original and modified total energy arises from the inconsistency between the original magnetostatic energy and the modified magnetostatic energy $r^2$.
%%%%%
%%%%%

\begin{figure}[!t]
\begin{center}
\subfigbottomskip = 1pt % set the distant between the figure in a row
\subfigcapskip = -2pt % set the distant between the subfigure and subtitile
\subfigure[$\Delta t = 10^{-13}$]{\label{Fig:SAVMRTUF-Diamond-1e13-Energy}
\includegraphics[width=69mm]{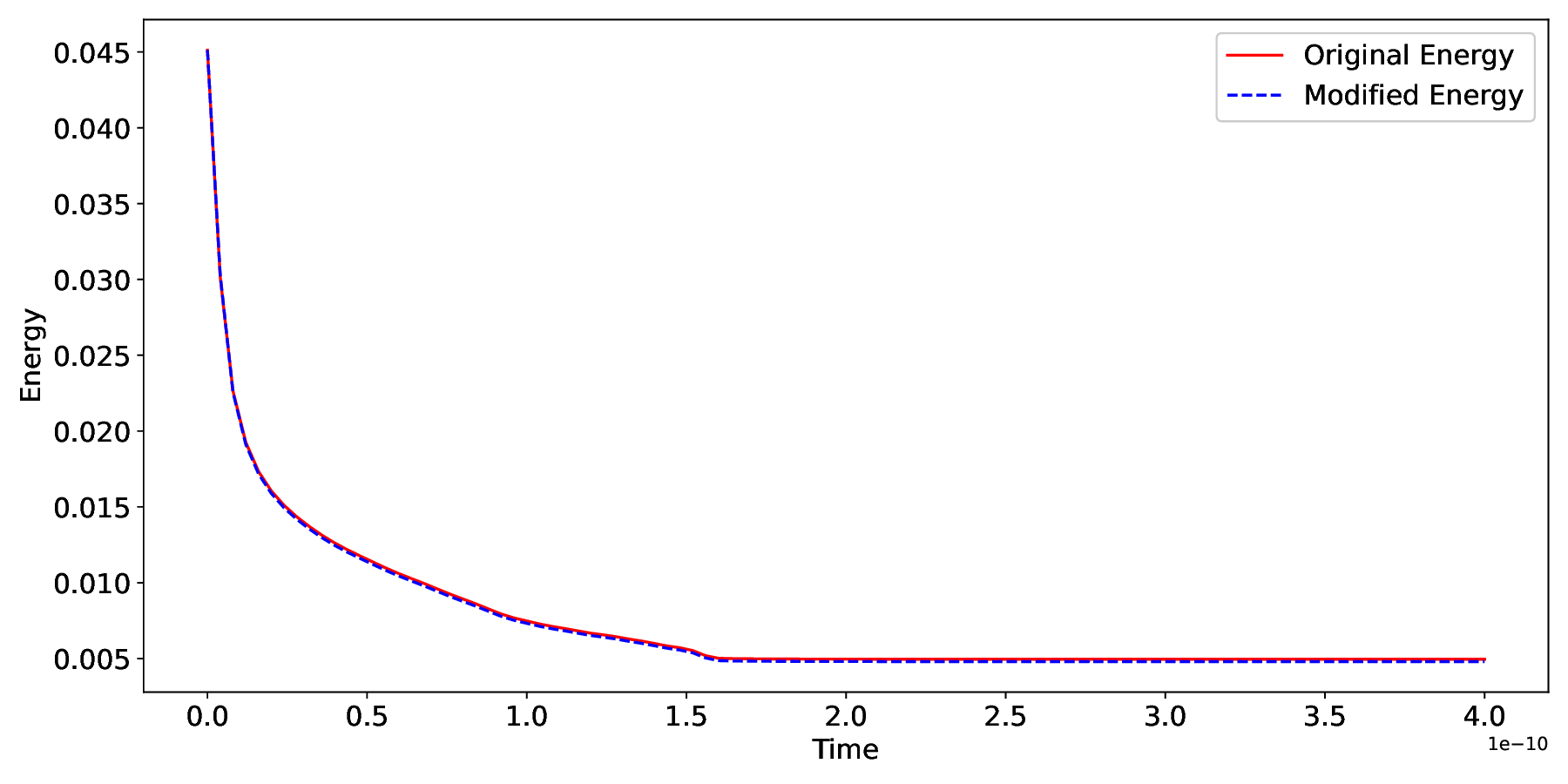}  %  fig 表示路径
}
\subfigure[$\Delta t = 10^{-14}$]{\label{Fig:SAVMRTUF-Diamond-1e14-Energy}
\includegraphics[width=69mm]{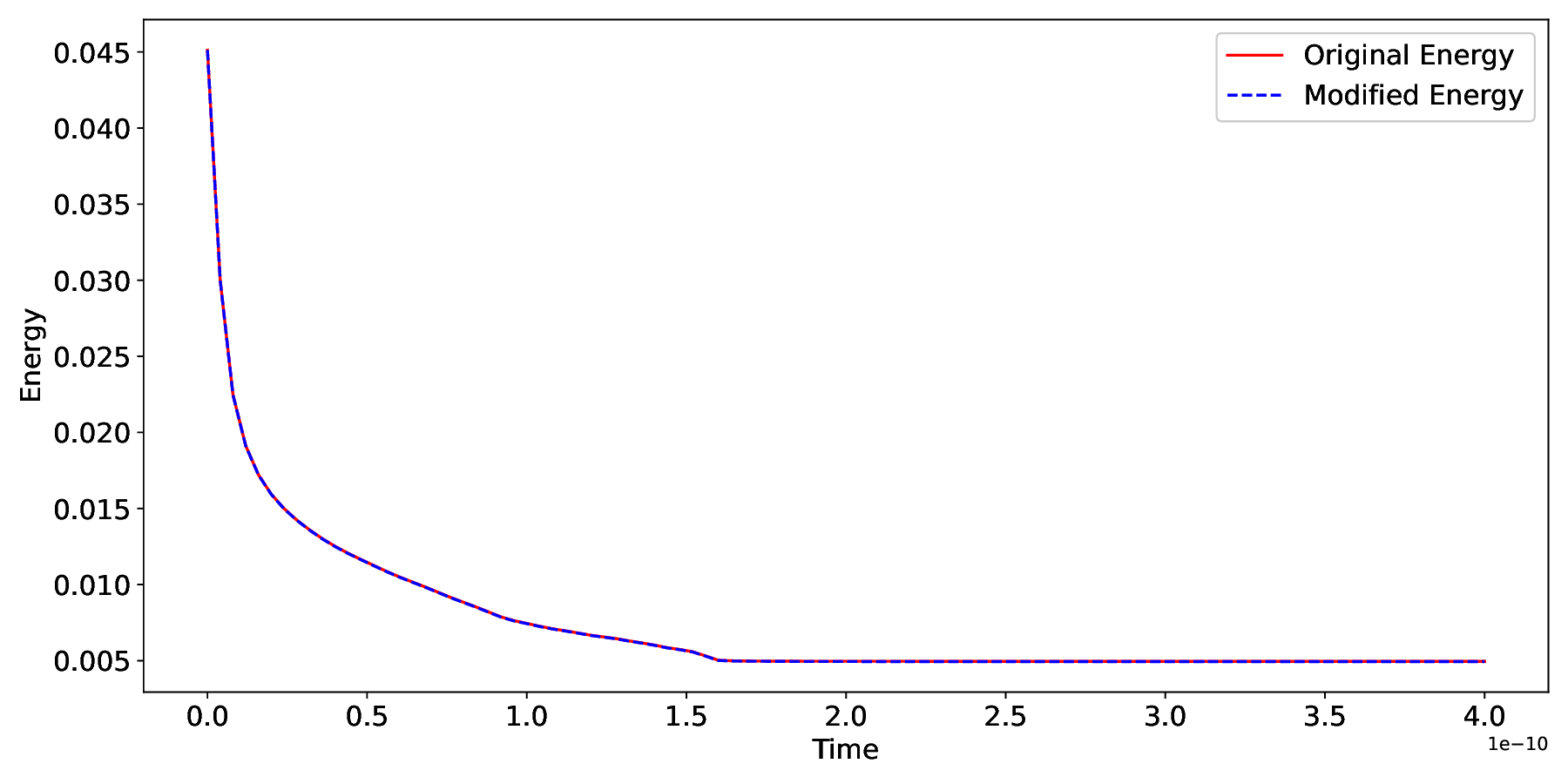}  %  fig 表示路径
}
\caption{Evolution of the original and modified energies over time under the SAV1 scheme with time steps of (a) $\Delta t = 10^{-13}$ and (b) $\Delta t = 10^{-14}$.}
\label{Fig:SAVMRTUF-Diamond-1e13}
\end{center}
\end{figure}

And then, two smaller temporal step sizes of $\Delta t = 10^{-13}$ and $10^{-14}$ are chosen.
Fig.~\ref{Fig:SAVMRTUF-Diamond-1e13} shows that the discrepancy between the original and modified energies vanishes with decreasing temporal step size.
%%%%%
%%%%%

\begin{figure}[!t]
\begin{center}
\subfigbottomskip = 1pt % set the distant between the figure in a row
\subfigcapskip = -2pt % set the distant between the subfigure and subtitile
\label{Fig:SAVRTUT-Diamond-AllEnergy}
\includegraphics[width=100mm]{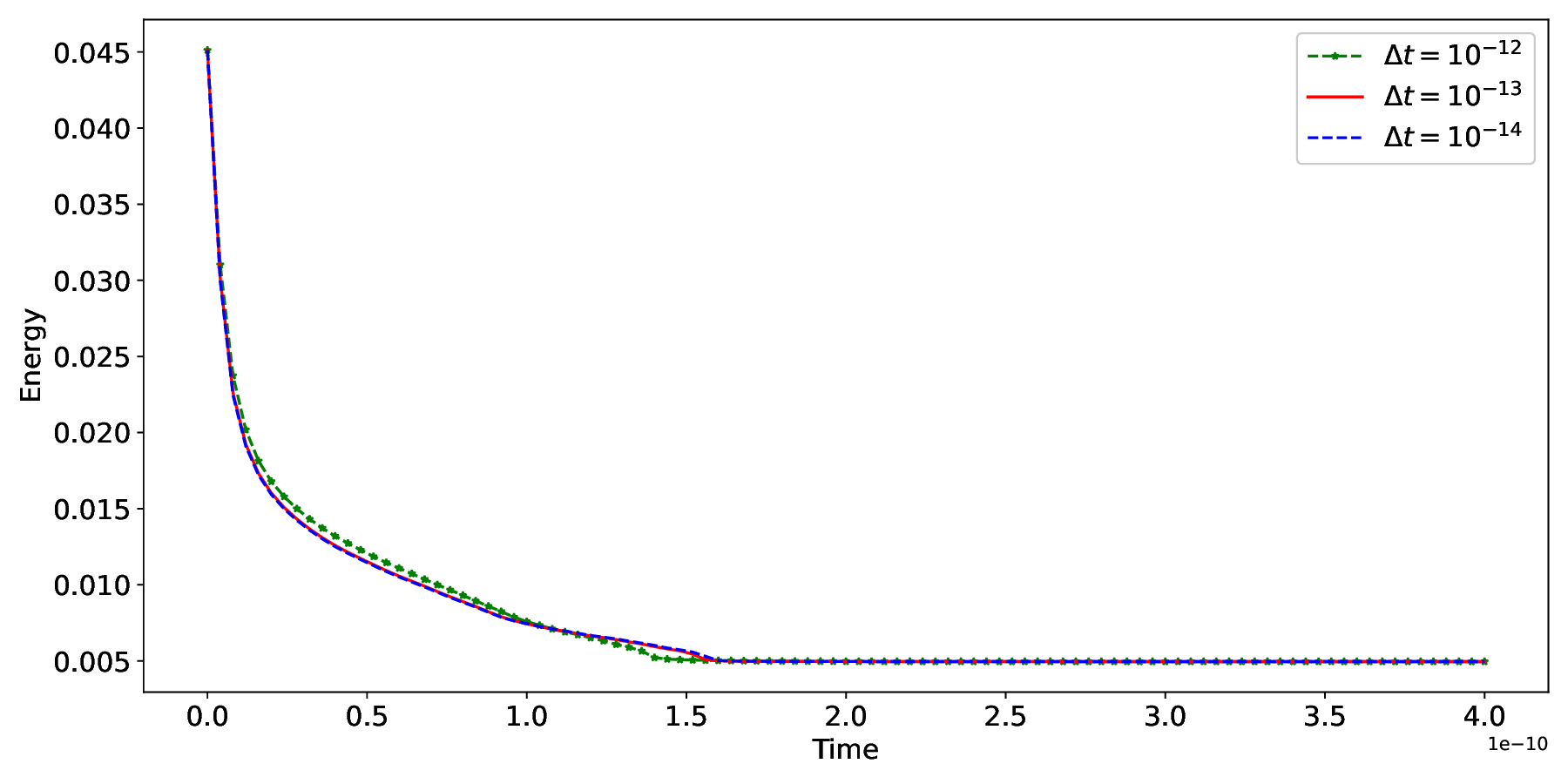}  %  fig 表示路径
\caption{Comparison of time-dependent total energy profiles for the SAV2 scheme with initial state \eqref{Eqn:IMDiamond} across a range of time steps ($\Delta t = 10^{-12}$, $10^{-13}$, and $10^{-14}$).}
\label{Fig:SAVRTUT-Diamond}
\end{center}
\end{figure}

Next, we test the effectiveness of SAV2 scheme \eqref{Eqn:RTUT} using Example \ref{Exa:1} with different temporal step sizes of $\Delta t = 10^{-12}$, $10^{-13}$, and $10^{-14}$.
In this scheme, the newly update equation \eqref{Eqn:RTUTd} for the scalar auxiliary variable $r$ enables that the original energy and the modified energy are equivalent.
From Fig.~\ref{Fig:SAVRTUT-Diamond}, it can be observed that all energies decrease smoothly and steadily, eventually reaching the equilibrium state.
Although the energy curve for the case of $\Delta t = 10^{-12}$ exhibits a small difference from the others, all of them converge to the same equilibrium. 
%%%%%
%%%%%

To further validate the effectiveness of both SAV1 and SAV2 schemes, we compare the corresponding steady-state energy with a reference value obtained by solving the micromagnetic dynamic equation --- Landau-Lifshitz-Gilbert (LLG) equation.
Precisely, the LLG equation is solved by using the following mid-point method with the parameters provided in Example \ref{Exa:1}:
\begin{equation}\label{Eqn:Mid-Point}
\dfrac{\h{m}^{n+1} - \h{m}^{n}}{\Delta t} = - \dfrac{ \h{m}^{n} + \h{m}^{n+1}}{2} \times \left( \dfrac{\h{h}_{\mathrm{eff}}^n + \h{h}_{\mathrm{eff}}^{n+1}}{2} - \alpha \dfrac{\h{m}^{n+1} - \h{m}^{n}}{\Delta t}  \right). 
\end{equation}
A sufficiently small temporal step size of $\Delta t = 10^{-14}$ is adopted.
Furthermore, we solve the nonlinear system in the mid-point method using explicit iteration with a tolerance error of $10^{-8}$.
After evolving for time $T = 2 \times 10^{-9}$, which is a sufficient time duration for the system to reach equilibrium, the resulting magnetization exhibits a diamond state as the steady state with a normalized energy of 0.004979. (In this section, all energies are normalized with respect to the stray field energy constant $K_d = \mu_0 M_s^2 /2$.)

The steady-state energies obtained by two SAV schemes and the mid-point method \eqref{Eqn:Mid-Point} are compared in Table \ref{Tab:MidPointSAV}.
In this table, all energies consistently demonstrate a high degree of similarity with the reference value, exhibiting a relative energy error of less than $0.5\%$. 
%%%
%%%
\begin{table}[!t]
	\centering\caption{The steady-state energies obtained by SAV1 and SAV2 schemes with the initial magnetization \eqref{Eqn:IMDiamond}, along with the corresponding relative error of energy.}
	\label{Tab:MidPointSAV}\vskip 0.1cm
	{\small
	\begin{tabular}{{ c | c | c |  c | c | c | c}}\hline
		Method   & \multicolumn{3}{c|}{SAV1} & \multicolumn{3}{c}{SAV2} \\ \hline 
$\Delta t$ & $ 10^{-12}$  &    $ 10^{-13}$   &  $ 10^{-14}$   & $ 10^{-12}$  &    $ 10^{-13}$   &  $ 10^{-14}$ \\   \hline
Energy & 0.004981 & 0.004955 & 0.004955 & 0.004957 & 0.004955 & 0.004955 \\ \hline
Relative error  & 0.04686\%  &  0.4747\%  & 0.4776\% & 0.4349\%  & 0.4778\%  & 0.4776\% \\
 \hline
	\end{tabular}
	}
\end{table}

\begin{figure}[!t]
\begin{center}
\subfigbottomskip = 1pt % set the distant between the figure in a row
\subfigcapskip = -2pt % set the distant between the subfigure and subtitile
\subfigure[Initial magnetization]{\label{Fig:SAVMRTUF-Initial-Landau}
\includegraphics[width=69mm]{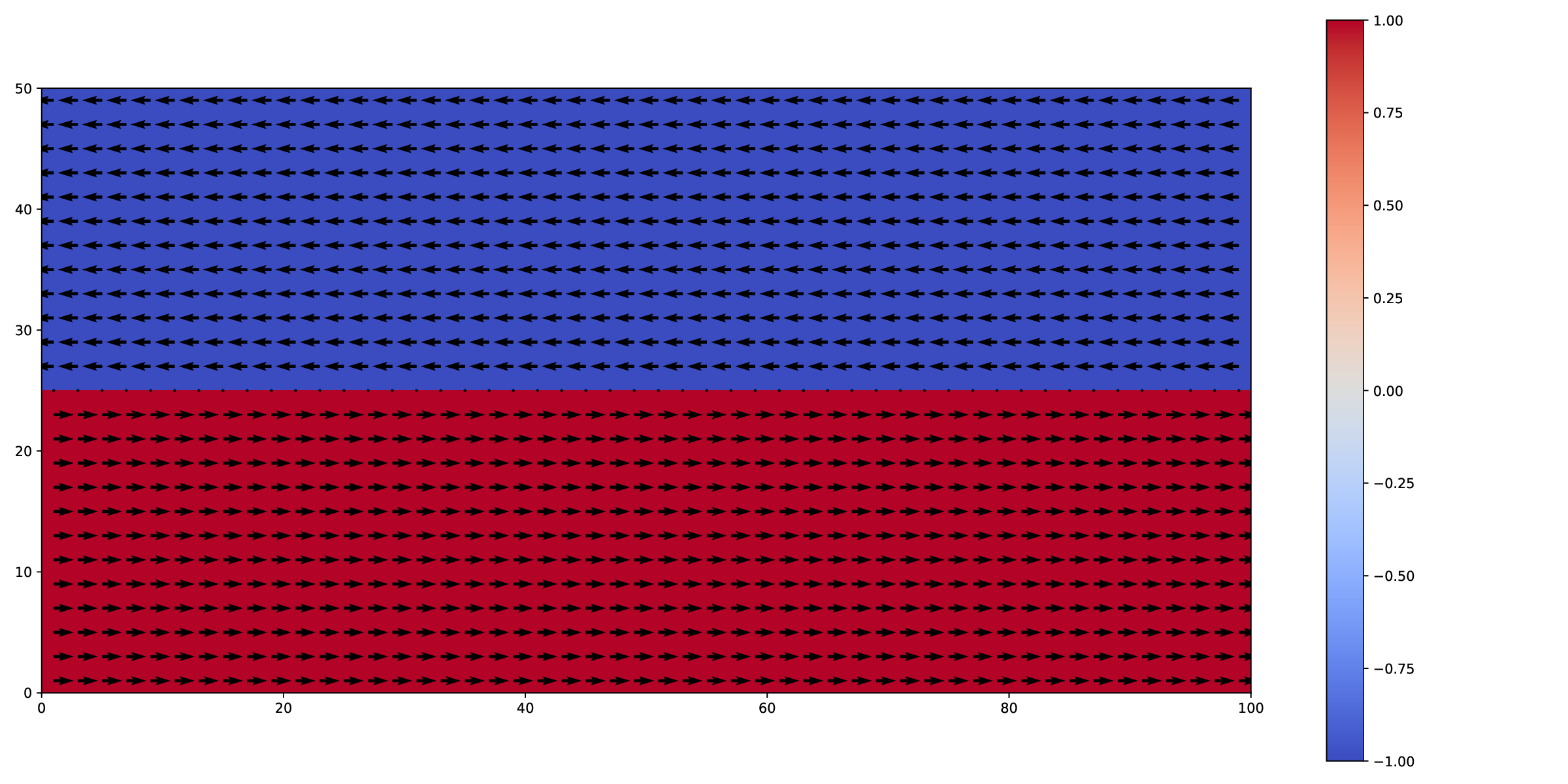}  %  fig 表示路径
}
\subfigure[Single cross-tie state]{\label{Fig:SAVMRTUF-SCT-M}
\includegraphics[width=69mm]{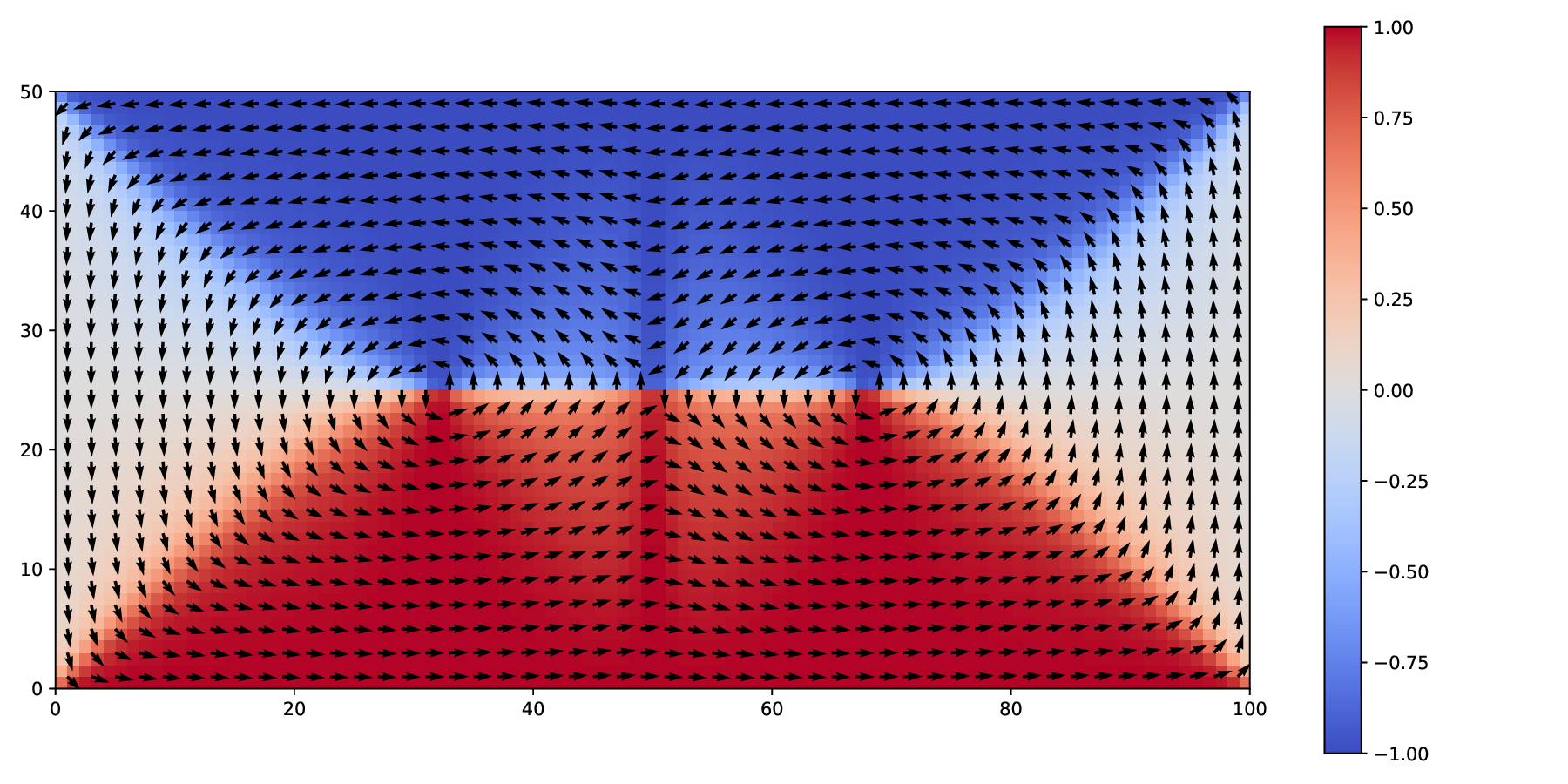}  %  fig 表示路径
}
\caption{(a) Initial magnetization \eqref{Eqn:IMSingleCrossTie}; (b) the corresponding steady-state named as single cross-tie state. The color mapping indicates the angle between the magnetization vector and the $x$-axis. }
\label{Fig:IMSCT}
\end{center}
\end{figure}

(ii)
To further validate the effectiveness of the SAV1 and SAV2 schemes, we conduct additional numerical simulations.
It is well known that different initial magnetizations can result in completely distinct steady states.
We reproduce the Example \ref{Exa:1}, but a distinct initial magnetization is considered at this time as
\begin{equation}\label{Eqn:IMSingleCrossTie}
\h{m}^0(x, y, z) = \left\{\begin{aligned}
& (1,0, 0), && \mbox{if $0  \leq y \leq 0.5 $}, \\
& (-1,0, 0), && \mbox{if $0.5  \leq y \leq 1 $}, 
\end{aligned}\right.
\end{equation}
which can lead to a single cross-tie state as the steady-state of magnetization, as shown in Fig.~\ref{Fig:IMSCT}. 
%%%

\begin{figure}[!th]
\begin{center}
\subfigbottomskip = 1pt % set the distant between the figure in a row
\subfigcapskip = -2pt % set the distant between the subfigure and subtitile
\includegraphics[width=100mm]{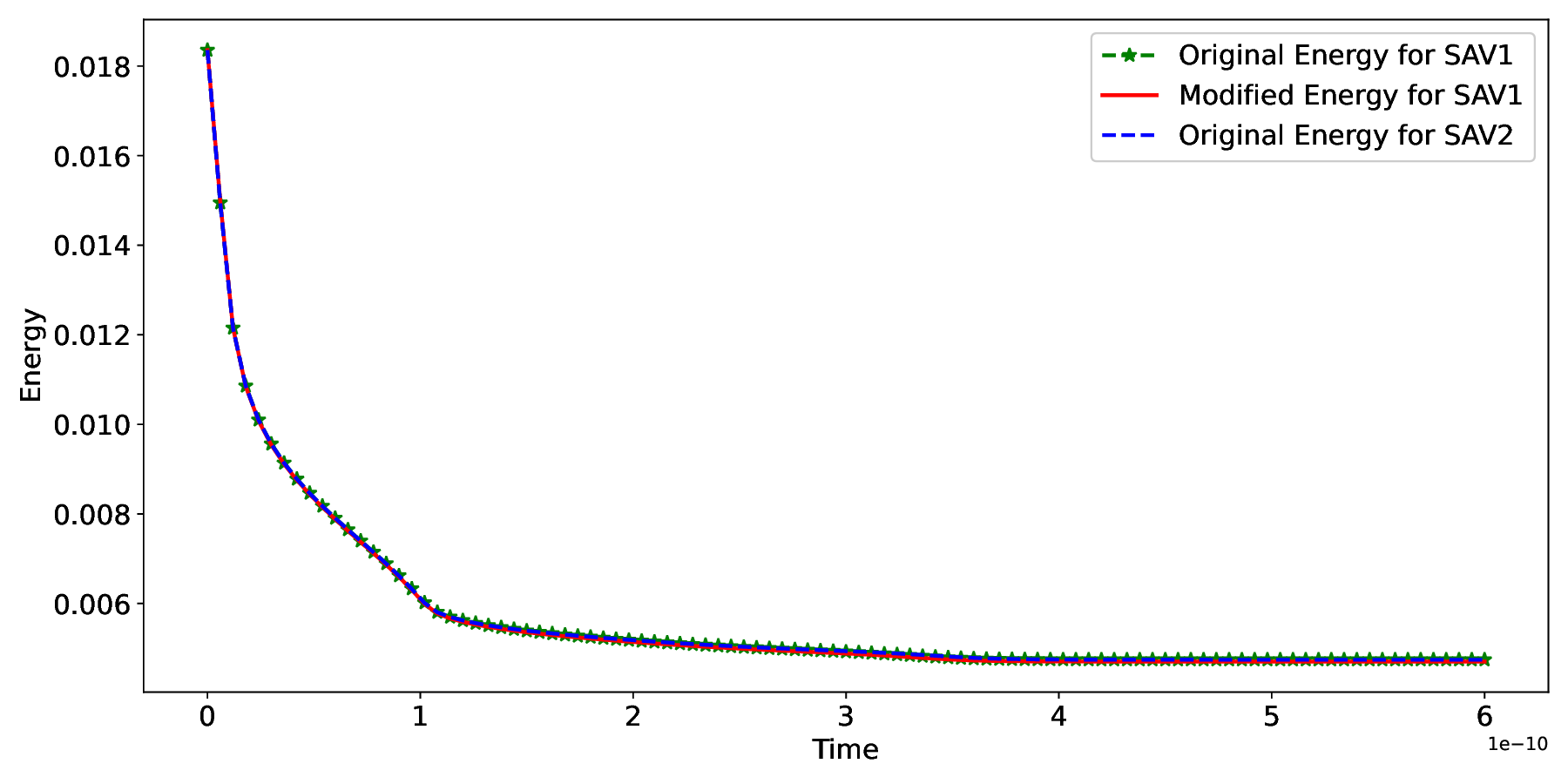}  %  fig 表示路径
\caption{Comparison of energies as functions of time by applying SAV1 and SAV2 schemes with initial magnetization \eqref{Eqn:IMSingleCrossTie} and $\Delta t = 10^{-13}$. }
\label{Fig:SAVMRTUF-SCT}
\end{center}
\end{figure}

We simulate the magnetization dynamics by applying SAV1 and SAV2 schemes with a temporal step size of $10^{-13}$ and simulation time $T = 6 \times 10^{-10}$.
It can be observed from Fig.~\ref{Fig:SAVMRTUF-SCT} that both schemes exhibit smooth and steady decrease in energy, and eventually reach same equilibrium.
Additionally, there is no obvious difference among these energy curves. 
The steady-state energies obtained by SAV1 and SAV2 schemes are compared in Table \ref{Tab:MidPointSAVSCT}.
It can be observed that all relative errors of energy, with the reference value of $0.004742$ obtained by the mid-point method \eqref{Eqn:Mid-Point}, are controlled within $0.03\%$.
%%%%
\begin{table}[!th]
	\centering\caption{The steady-state energies obtained by SAV1 and SAV2 schemes with the initial magnetization \eqref{Eqn:IMSingleCrossTie} and \eqref{Eqn:IMDoubleCrossTie}, and $\Delta t = 10^{-13}$, along with the corresponding relative error of energy.}
	\label{Tab:MidPointSAVSCT}\vskip 0.1cm
	\begin{tabular}{{ c | c | c |  c | c }}\hline
Steady-state   &   \multicolumn{2}{c|}{Single Cross-tie} & \multicolumn{2}{c}{Double Cross-tie} \\  \hline
		Methods   & SAV1 & SAV2 & SAV1 & SAV2 \\ \hline 
Energy &   0.004743 & 0.004743 & 0.005068  & 0.005067 \\  \hline
Relative error & 0.02557\% & 0.02509\% & 0.9587\%   & 0.9280\% \\  
 \hline
	\end{tabular}
\end{table}

\begin{figure}[!ht]
\begin{center}
\subfigbottomskip = 1pt % set the distant between the figure in a row
\subfigcapskip = -2pt % set the distant between the subfigure and subtitile
\subfigure[Initial magnetization]{\label{Fig:SAVMRTUF-Initial-DCT}
\includegraphics[width=69mm]{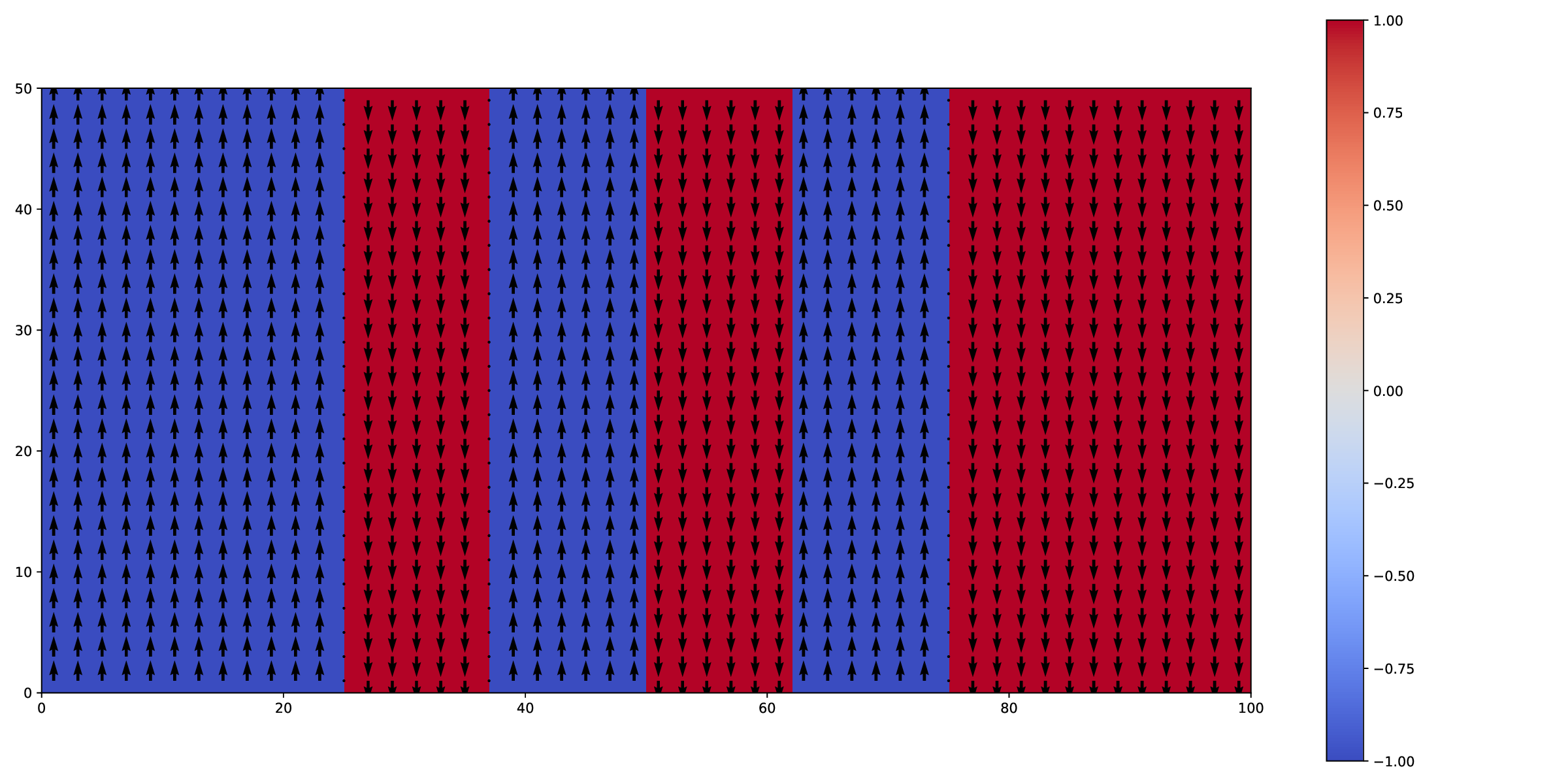}  %  fig 表示路径
}
\subfigure[Double cross-tie state]{\label{Fig:SAVMRTUF-DCT-M}
\includegraphics[width=69mm]{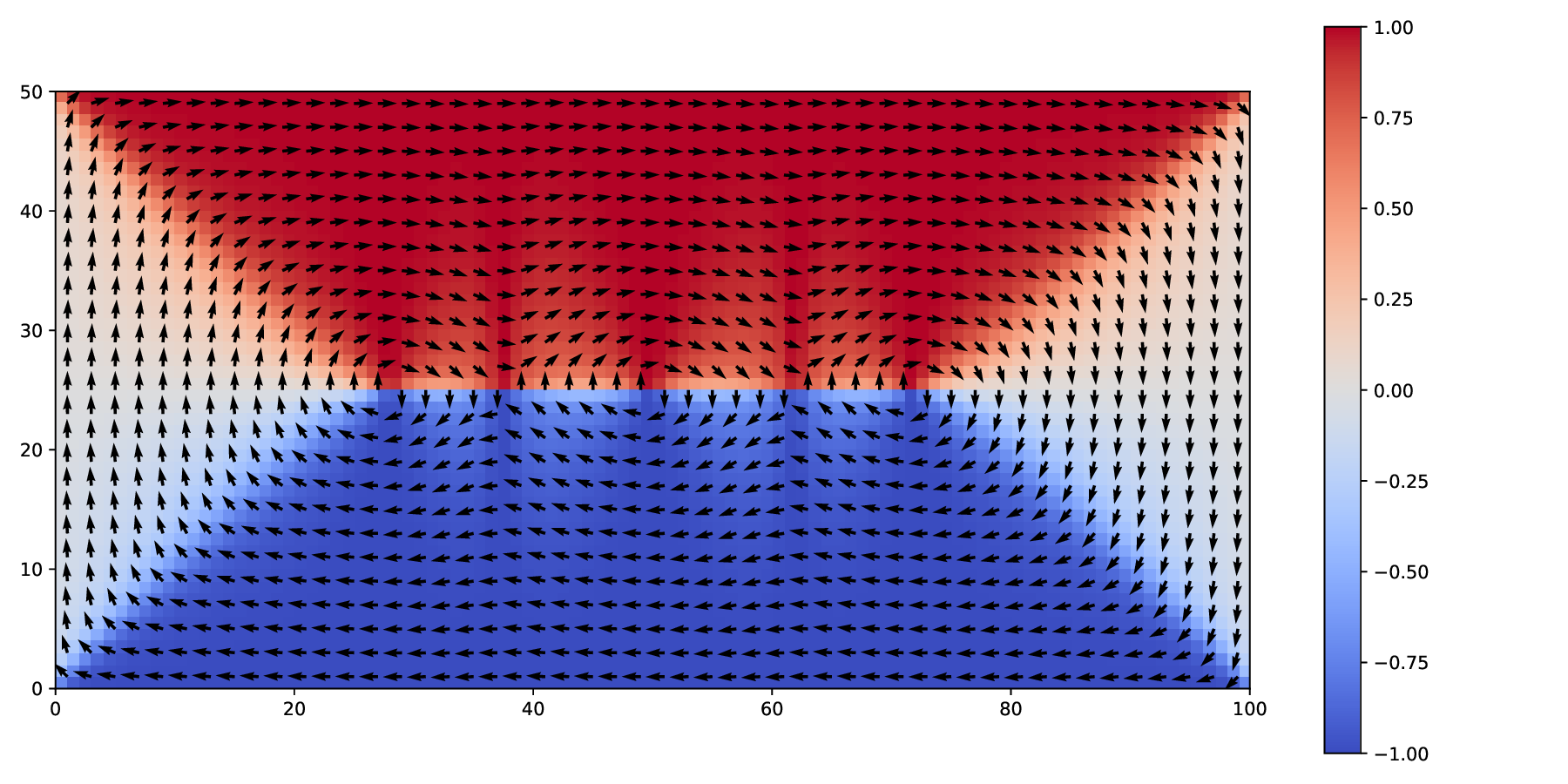}  %  fig 表示路径
}
\caption{(a) Initial magnetization \eqref{Eqn:IMDoubleCrossTie}; (b) the corresponding steady-state named as double cross-tie state. The background color in Figure (a) represents the angle between the magnetization and the $y$-axis, while in Figure (b), it indicates the angle between the magnetization vector and the $x$-axis. }
\label{Fig:IMDCT}
\end{center}
\end{figure}

\begin{figure}[!t]
\begin{center}
\subfigbottomskip = 1pt % set the distant between the figure in a row
\subfigcapskip = -2pt % set the distant between the subfigure and subtitile
\includegraphics[width=100mm]{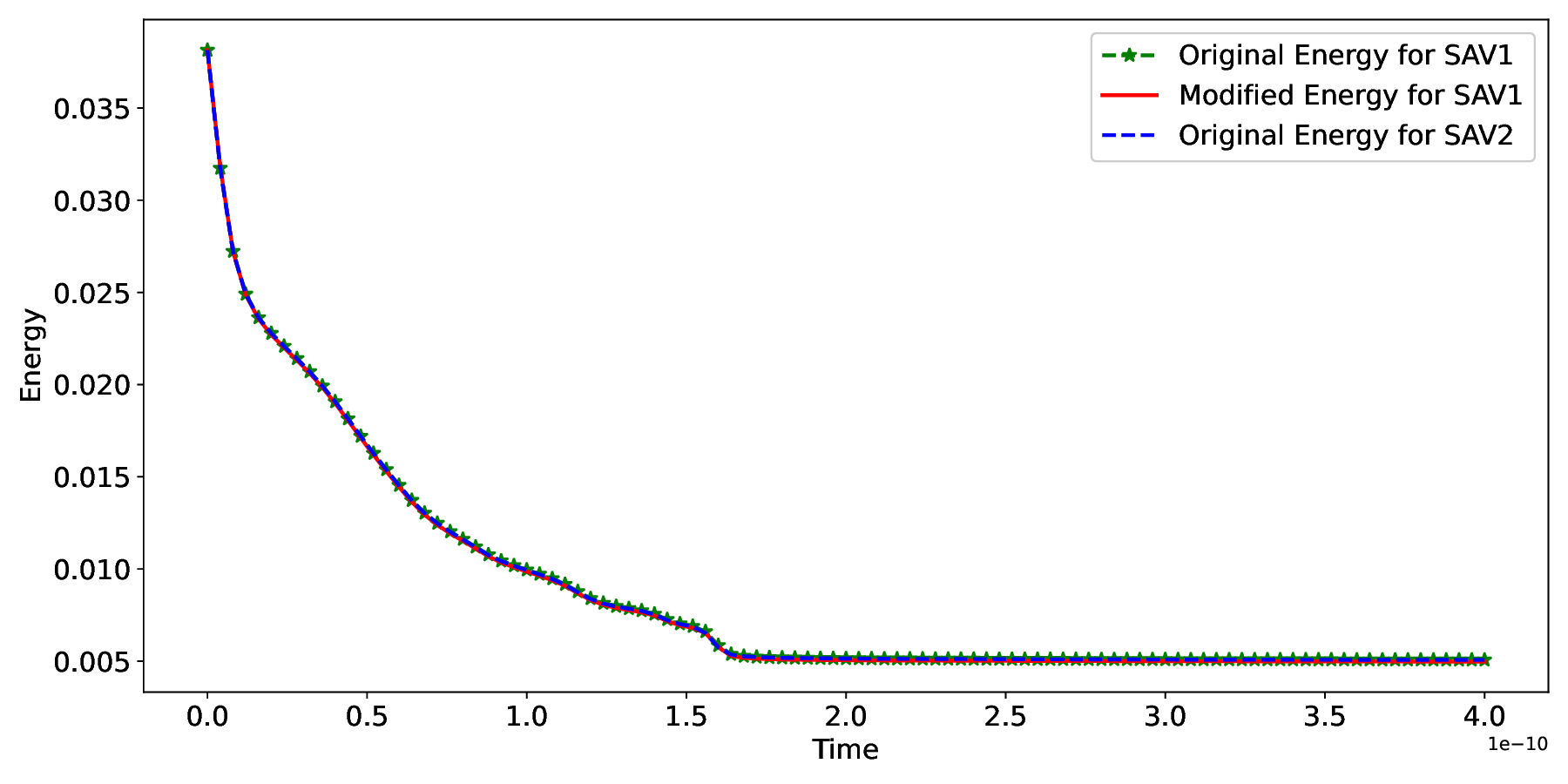}  %  fig 表示路径
\caption{Comparison of energies as functions of time by applying SAV1 and SAV2 schemes with initial magnetization \eqref{Eqn:IMDoubleCrossTie} and $\Delta t = 10^{-13}$. }
\label{Fig:SAVMRTUF-DCT}
\end{center}
\end{figure}

(iii) At last, we simulate the Example \ref{Exa:1} with another initial magnetization
\begin{equation}\label{Eqn:IMDoubleCrossTie}
\h{m}^0(x, y, z) = \left\{\begin{aligned}
& (0,1, 0), && \mbox{if $x \in (0, 0.5) \cup (0.75, 1) \cup (1.25 , 1.5)$}, \\
& (0,-1, 0), && \mbox{if $x \in (0.5, 0.75) \cup (1 , 1.25) \cup (1.5 , 2)$}, 
\end{aligned}\right.
\end{equation}
to achieve a double cross-tie state, as shown in Fig.~\ref{Fig:IMDCT}. 
%%%
Similarly, it can be observed from Fig.~\ref{Fig:SAVMRTUF-DCT} that all energy curves obtained by SAV1 and SAV2 schemes steadily decrease, eventually reaching equilibrium, with no significant discrepancies among them.
Moreover, as observed in Table \ref{Tab:MidPointSAVSCT}, the relative errors of energy, with the reference value of $0.004742$ obtained by the mid-point method \eqref{Eqn:Mid-Point}, are consistently controlled within $1\%$.
%%%%%
%%%%

From the above numerical simulations, it can be seen that all the energies consistently decrease and eventually reach an equilibrium.
Although the projection of $\h{m}^*$ and the change in the update of the scalar auxiliary variable bring challenges in theoretically presenting the energy dissipation analysis of both SAV schemes (See Remark \ref{Rem:Energy}),
throughout the magnetization evolution, numerical tests confirm that both SAV schemes maintain consistent energy dissipation.

(iv) At last, some convergence test for both SAV1 and SAV2 schemes are performed by Example \ref{Exa:1}.
The energies of the steady-state obtained by these two schemes are investigated by varying the mesh sizes and temporal step sizes.
To ensure the attainment of steady state, the numerical simulations are implemented with a simulation time duration $T = 4\times 10^{-10}$.
This value is chosen, ensuring adequate time for the system to attain equilibrium.

\begin{figure}[!t]
\begin{center}
\subfigbottomskip = 1pt % set the distant between the figure in a row
\subfigcapskip = -2pt % set the distant between the subfigure and subtitile
\subfigure[Iteration numbers VS Energy]{\label{Fig:MRIterationNumbersVSEnergy}
\includegraphics[width=69mm]{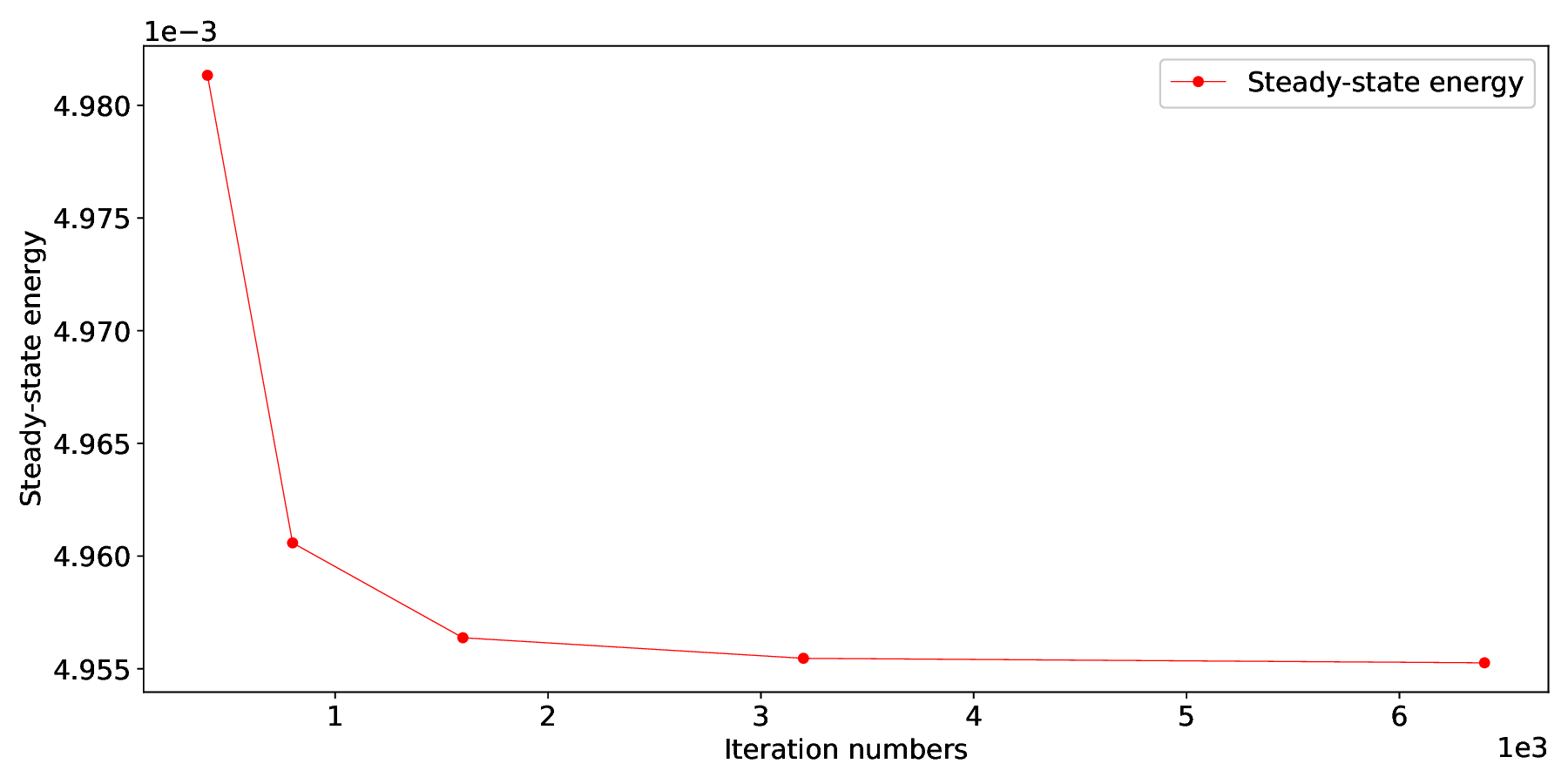}  %  fig 表示路径
}
\subfigure[Number of degrees of freedom VS Energy]{\label{Fig:MRNumbersOfFreedomVSEnergy}
\includegraphics[width=69mm]{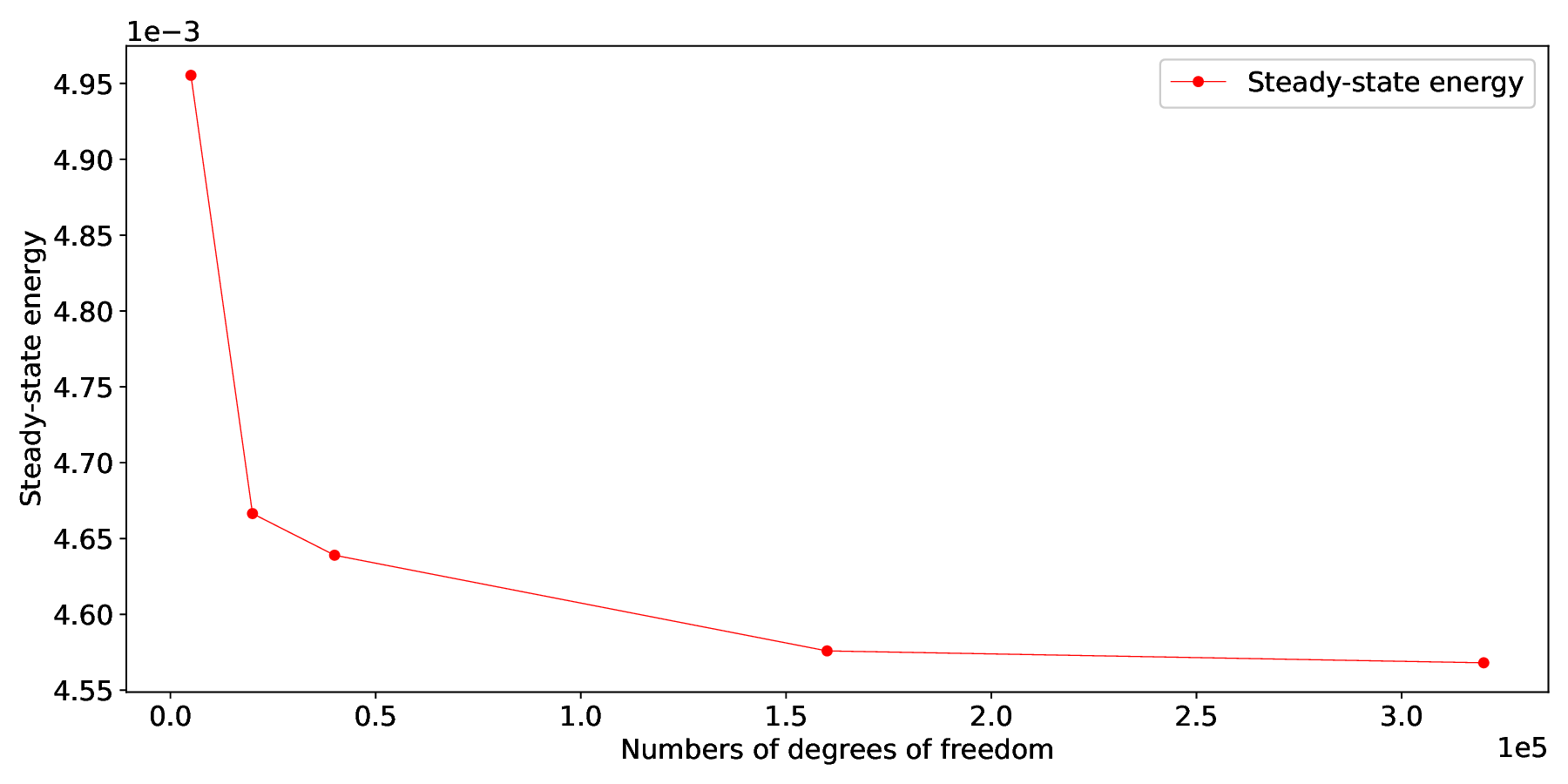}  %  fig 表示路径
}
\caption{Convergence test for SAV1 scheme. Left: energies of steady-state with different iteration numbers. 
Right: energies of steady-state with different mumber of degrees of freedom.}
\label{Fig:MRConvergenceTest}
\end{center}
\end{figure}

\begin{figure}[!t]
\begin{center}
\subfigbottomskip = 1pt % set the distant between the figure in a row
\subfigcapskip = -2pt % set the distant between the subfigure and subtitile
\subfigure[Iteration numbers VS Energy]{\label{Fig:IterationNumbersVSEnergy}
\includegraphics[width=69mm]{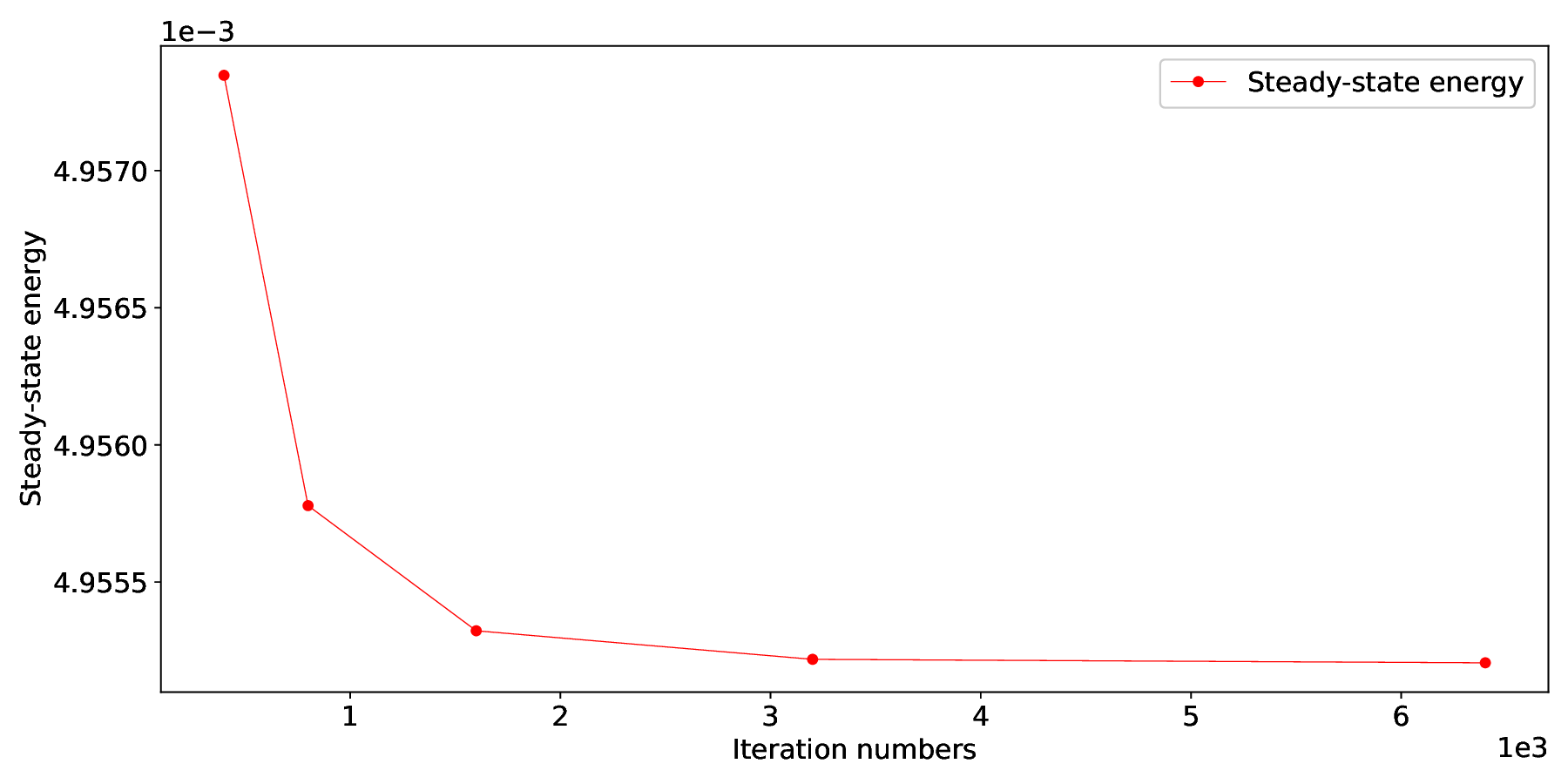}  %  fig 表示路径
}
\subfigure[Number of degrees of freedom VS Energy]{\label{Fig:NumbersOfFreedomVSEnergy}
\includegraphics[width=69mm]{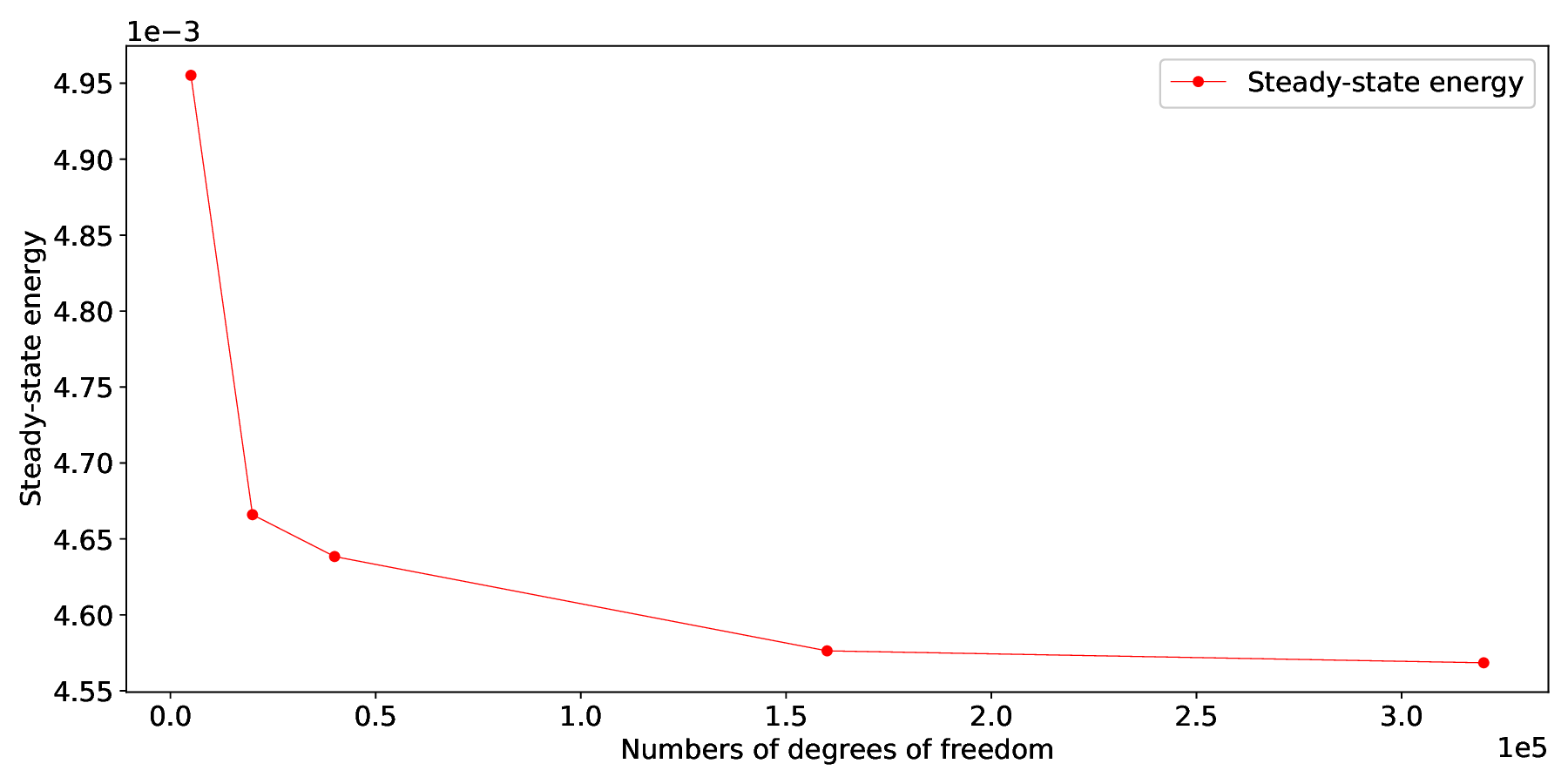}  %  fig 表示路径
}
\caption{Convergence test for SAV2. Left: energies of steady-state with different iteration numbers. 
Right: energies of steady-state with different number of degrees of freedom.}
\label{Fig:ConvergenceTest}
\end{center}
\end{figure}

Firstly, we fix the spacial mesh to $(100, 50, 1)$ and investigate the energies of steady-state obtained by  SAV1 and SAV2 schemes with different temporal step sizes.
From Figs.~\ref{Fig:MRIterationNumbersVSEnergy} and \ref{Fig:IterationNumbersVSEnergy}, we observe a gradual convergence of the steady-state energies as the iteration numbers increase (which means the temporal step size decreases).
%%%%
%%%%%

And then, the temporal step size is fixed as $10^{-13}$, and different numbers of degrees of freedom are adopted.
From Figs.~\ref{Fig:MRNumbersOfFreedomVSEnergy} and \ref{Fig:NumbersOfFreedomVSEnergy}, it also can be observed that the steady-state energies gradually stabilize as the number of degrees of freedom increases.

\subsection{Efficiency}\label{Sec:Efficiency}

In this part, the efficiency of the SAV1 and SAV2 schemes is validated by comparing them with commonly employed numerical methods for solving the flow equation \eqref{Eqn:GF} and the LLG equation.

(i) We compare the efficiency of the FEP scheme \eqref{Eqn:FE}, the BEP scheme \eqref{Eqn:BE}, and the SAV schemes \eqref{Eqn:SAVTD} and \eqref{Eqn:RTUT} for solving the flow equation \eqref{Eqn:GF}.
Example \ref{Exa:1} is implemented with the initial magnetization \eqref{Eqn:IMDiamond} to achieve the diamond state.
Similar to the previous simulation, a uniform spacial mesh of $(100, 50, 1)$ and total simulation time of $T = 4 \times 10^{-10}$ are employed.
An explicit iterative solution was implemented for the BEP scheme with a tolerance error of $10^{-8}$.
The CPU time is exhibited only in the case where the relative error of energy, with the reference value of $0.004979$, is less than $1\%$.

We present the CPU times of the FEP, BEP, and SAV schemes with different temporal step sizes in Table \ref{Tab:CPUGF}.
It can be seen that the BEP scheme suffers from tedious computation time since it needs to be solved using an iterative method as demonstrated in Section \ref{Sec:BCFE}.
Additionally, the CPU times of SAV1, SAV2 and FEP schemes are comparable.
Although the time cost of SAV1 and SAV2 schemes are slightly higher than that of the FEP scheme, SAV schemes allow for more flexible temporal step sizes.
The SAV schemes still perform well when the temporal step size is $1.42 \times 10^{-12}$, whereas the FEP scheme fails. 
In conclusion, both proposed SAV schemes are not only efficient but also have fewer limitations for the temporal step size.
%%%%
\begin{table}[htp]
	\centering\caption{Comparison of CPU times (s) for different methods to achieve the diamond state by solving flow equation \eqref{Eqn:GF} with different $\Delta t$.}
	\label{Tab:CPUGF}
	\vspace{0.2cm}
	\begin{tabular}{{ c | c | c | c | c | c|  c  }}\hline
		\diagbox{\footnotesize Scheme}{\footnotesize $\Delta t$}  & $1.42 \times 10^{-12}$ & $10^{-12} $  & $ 5 \times 10^{-13}$  & $ 10^{-13}$ & $  5 \times 10^{-14}$  & $  10^{-14}$ \\ \hline
				FEP & - &  - & 7.32 &   35.56 & 64.25  &  315.24 \\
		BEP & - & -  & 154.10 &  264.03 &  424.36  &   1356.14 \\
		{\bf SAV1	} & {\bf 3.39} & {\bf 5.16} &  {\bf  9.15} & {\bf 47.97} &  {\bf 99.96} &  {\bf 499.04}  \\
		{\bf SAV2	} & {\bf 3.49}& {\bf 5.08} &  {\bf 10.62} & {\bf 49.15} &  {\bf  99.84} &  {\bf 508.78}  \\  
 \hline
	\end{tabular}
\end{table}

(ii) To further show the high-efficiency of SAV schemes, we compare SAV1 and SAV2 schemes with the mid-point method \eqref{Eqn:Mid-Point}, the backward Euler method and the forward Euler method of LLG equation.
The backward Euler method and forward Euler method of LLG equation are described as follows:
\begin{eqnarray}\label{Eqn:LLGBE}
& \dfrac{\h{m}^{n+1} - \h{m}^{n}}{\Delta t} = -\h{m}^{n+1} \times \left( \h{h}_{\mathrm{eff}}^{n+1} - \alpha \dfrac{\h{m}^{n+1} - \h{m}^{n}}{\Delta t}  \right),
\\ \label{Eqn:LLGFE}
& \dfrac{\h{m}^{n+1} - \h{m}^{n}}{\Delta t} = -\h{m}^{n} \times \left( \h{h}_{\mathrm{eff}}^{n} - \alpha \dfrac{\h{m}^{n+1} - \h{m}^{n}}{\Delta t}  \right). 
\end{eqnarray}

The CPU time tests for the LLG equation are conducted over a duration of time $T = 2\times 10^{-9}$ to ensure the attainment of a steady state. 
The optimal CPU times of the mid-point method \eqref{Eqn:Mid-Point}, the backward Euler method \eqref{Eqn:LLGBE}, and the forward Euler method \eqref{Eqn:LLGFE} are presented in Table \ref{Tab:CPUGFLLG}.
It is evident that the CPU times of the SAV schemes \eqref{Eqn:SAVTD} and \eqref{Eqn:RTUT} are significantly lower than those of solving the midpoint, backward Euler, and forward Euler methods of the LLG equation.
This is because in the midpoint and backward Euler methods, nonlinear systems need to be solved, while the forward Euler method requires solving a linear system with variable coefficient.
However, SAV1 and SAV2 schemes only require solving two linear systems with constant coefficient, which can be efficiently solved by using the Discrete Cosine Transform.
Furthermore, the simulation time required to reach the steady-state for the LLG equation is significantly longer than that for the Gibbs energy minimization problem.
Therefore, both SAV1 and SAV2 schemes have a significant advantage over solving the LLG equation in terms of efficiency for micromagnetic ground state problems.
%%%%
%%%%
\begin{table}[htp]
	\centering\caption{Comparison of the optimal CPU times (s) for mid-point method \eqref{Eqn:Mid-Point}, backward Euler method \eqref{Eqn:LLGBE}, and forward Euler method \eqref{Eqn:LLGFE} to achieve the diamond state.}
	\label{Tab:CPUGFLLG}
	\vspace{0.2cm}
	{\footnotesize
	\begin{tabular}{{ c | c | c |  c  }}\hline
		Methods   & mid-point \eqref{Eqn:Mid-Point}  & backward Euler \eqref{Eqn:LLGBE} &  forward Euler \eqref{Eqn:LLGFE} \\ \hline
	CPU time	   & 88.49 ($\Delta t = 4 \times 10^{-12}$) & 639.69 ($\Delta t = 1.25\times 10^{-13}$)& 162.04 ($\Delta t = 1.25\times 10^{-13}$) \\ \hline
	\end{tabular}
	}
\end{table}

\section{Conclusion}\label{Sec:Conclusion}

In this paper, two efficient SAV schemes are proposed to solve the flow equation \eqref{Eqn:GF} for Gibbs energy minimization, aiming to achieve the ground state of the micromagnetic structure.
The SAV1 scheme \eqref{Eqn:SAVTD} is proposed for the flow equation \eqref{Eqn:GF} by combining the original SAV method with a projection method that ensures the constant length of magnetization.
To enhance the correspondence between the scalar auxiliary variable and the actual energy, we modify the update approach of the scalar auxiliary variable and propose SAV2 scheme \eqref{Eqn:RTUT}.
Both SAV schemes are demonstrated to partially preserve energy dissipation and possess high computational efficiency.
The primary computational expense lies in solving two linear systems that possess constant coefficients.
Additionally, the utilization of the Discrete Cosine Transform can significantly enhance the efficiency of solving these two linear systems.
The numerical experiments reveal that both SAV1 and SAV2 schemes have a significant advantage over other numerical methods in terms of efficiency.


\vspace{0.6cm}


\begin{thebibliography}{10}

\bibitem{AkrivisGLiBY19:A3703}
G.~Akrivis, B.~Li, and D.~Li.
\newblock Energy-decaying extrapolated {RK}–{SAV} methods for the
  {Allen}–{Cahn} and {Cahn}–{Hilliard} equations.
\newblock {\em SIAM J. Sci. Comput.}, 41(6):A3703--A3727, 2019.

\bibitem{BernadouMDepeyreS02:1599}
M.~Bernadou, S.~Depeyre, and S.~He.
\newblock Numerical modelization of magnetic materials based on energy
  minimization.
\newblock {\em AIP Conf. Proc.}, 615(1):1599--1606, 2002.

\bibitem{BernadouMDepeyreS02:1018}
M.~Bernadou, S.~Depeyre, S.~He, and P.~Meilland.
\newblock Numerical simulation of magnetic microstructure based on energy
  minimization.
\newblock {\em J. Magn. Magn. Mater.}, 242–245:1018--1020, 2002.

\bibitem{BernadouMDepeyreS03:86}
M.~Bernadou, S.~Depeyre, S.~He, and P.~Meilland.
\newblock A numerical simulation of ferromagnetic material at microscopic
  scale.
\newblock In {\em 13th International Conference on Adaptive Structures and
  Technologies}, pages 86--93. CRC Press, 2003.

\bibitem{BernadouMDepeyreS04:164}
M.~Bernadou, S.~Depeyre, S.~He, and P.~Meilland.
\newblock Numerical approximation of ferromagnetic materials.
\newblock In {\em Smart Structures and Materials 2004: Modeling, Signal
  Processing, and Control}, volume 5383, pages 164--173. SPIE, 2004.

\bibitem{BernadouMHeS98:512}
M.~Bernadou and S.~He.
\newblock On the computation of hysteresis and closure domains in
  micromagnetism.
\newblock In {\em Smart Structures and Materials 1998: Mathematics and Control
  in Smart Structures}, volume 3323, pages 512--519. SPIE, 1998.

\bibitem{BernadouMHeS}
M.~Bernadou and S.~He.
\newblock Numerical approximation of unstressed or prestressed magnetostrictive
  materials.
\newblock In {\em Smart Structures and Materials 1999: Mathematics and Control
  in Smart Structures}, volume 3667, pages 101--109. SPIE, 1999.

\bibitem{BlueJLScheinfeinMR91:4778}
J.~L. Blue and M.~R. Scheinfein.
\newblock Using multipoles decreases computation time for magnetostatic
  self-energy.
\newblock {\em IEEE Trans. Magn.}, 27(6):4778--4780, 1991.

\bibitem{CarstensenCPraetoriusD05:2633}
C.~Carstensen and D.~Praetorius.
\newblock Numerical analysis for a macroscopic model in micromagnetics.
\newblock {\em SIAM J. Numer. Anal.}, 42(6):2633--2651, 2005.

\bibitem{CarstensenCProhlA01:65}
C.~Carstensen and A.~Prohl.
\newblock Numerical analysis of relaxed micromagnetics by penalised finite
  elements.
\newblock {\em Numer. Math.}, 90:65--99, 2001.

\bibitem{ChenCYangX19:35}
C.~Chen and X.~Yang.
\newblock Fast, provably unconditionally energy stable, and second-order
  accurate algorithms for the anisotropic {Cahn}–{Hilliard} model.
\newblock {\em Comput. Methods Appl. Mech. Eng.}, 351:35--59, 2019.

\bibitem{ChenJRWangC21:55}
J.~Chen, C.~Wang, and C.~Xie.
\newblock Convergence analysis of a second-order semi-implicit projection
  method for {Landau}-{Lifshitz} equation.
\newblock {\em Appl. Numer. Math.}, 168:55--74, 2021.

\bibitem{ChengQLiuC21:113532}
Q.~Cheng, C.~Liu, and J.~Shen.
\newblock Generalized {SAV} approaches for gradient systems.
\newblock {\em J. Comput. Appl. Math.}, 394:113532, 2021.

\bibitem{ChengQShenJ18:A3982}
Q.~Cheng and J.~Shen.
\newblock Multiple scalar auxiliary variable ({MSAV}) approach and its
  application to the phase-field vesicle membrane model.
\newblock {\em SIAM J. Sci. Comput.}, 40(6):A3982--A4006, 2018.

\bibitem{CuiZYangL24:}
Z.~Cui, L.~Yang, J.~Wu, and G.~Hu.
\newblock A treecode algorithm for the {Poisson} equation in a general domain
  with unstructured grids.
\newblock {\em Numer. Algorithms}, 2024.
\newblock Publish Online: https://doi.org/10.1007/s11075-024-01888-8.

\bibitem{EWNWangXP01:1647}
W.~E and X.~Wang.
\newblock Numerical methods for the {Landau}-{Lifshitz} equation.
\newblock {\em SIAM J. Numer. Anal.}, 38(5):1647--1665, 2001.

\bibitem{ExlLAuzingerW12:2840}
L.~Exl, W.~Auzinger, S.~Bance, M.~Gusenbauer, F.~Reichel, and T.~Schrefl.
\newblock Fast stray field computation on tensor grids.
\newblock {\em J. Comput. Phys.}, 231(7):2840--2850, 2012.

\bibitem{ExlLBanceS14:17D118}
L.~Exl, S.~Bance, F.~Reichel, T.~Schrefl, H.~P. Stimming, and N.~J. Mauser.
\newblock {LaBonte's} method revisited: An effective steepest descent method
  for micromagnetic energy minimization.
\newblock {\em J. Appl. Phys.}, 115:17D118, 2014.

\bibitem{ExlLFischbacherJ19:179}
L.~Exl, J.~Fischbacher, A.~Kovacs, H.~Oezelt, M.~Gusenbauer, and T.~Schrefl.
\newblock Preconditioned nonlinear conjugate gradient method for micromagnetic
  energy minimization.
\newblock {\em Comput. Phys. Commun.}, 235:179--186, 2019.

\bibitem{FischbacherJKovacsA17:045310}
J.~Fischbacher, A.~Kovacs, H.~Oezelt, T.~Schrefl, L.~Exl, J.~Fidler, D.~Suess,
  N.~Sakuma, M.~Yano, A.~Kato, T.~Shoji, and A.~Manabe.
\newblock Nonlinear conjugate gradient methods in micromagnetics.
\newblock {\em AIP Adv.}, 7:045310, 2017.

\bibitem{FuruyaAFujisakiJ15:1}
A.~Furuya, J.~Fujisaki, K.~Shimizu, Y.~Uehara, T.~Ataka, T.~Tanaka, and
  H.~Oshima.
\newblock Semi-implicit steepest descent method for energy minimization and its
  application to micromagnetic simulation of permanent magnets.
\newblock {\em IEEE Trans. Magn.}, 51(11):1--4, 2015.

\bibitem{GarciaCervera2007:Review}
C.~J. Garc\'{ı}a-Cervera.
\newblock Numerical micromagnetics: a review.
\newblock {\em Bol. Soc. Esp. Mat. Apl.}, 39:103--135, 2007.

\bibitem{HeJYangL24:1179}
J.~He, L.~Yang, and J.~Zhan.
\newblock Temporal high-order accurate numerical scheme for the
  {Landau}–{Lifshitz}–{Gilbert} equation.
\newblock {\em Mathematics}, 12(8):1179, 2024.

\bibitem{HeJYang}
J.~He, L.~Yang, J.~Chen, and J.~Zhan.
\newblock Structure-preserving numerical schemes for the ground state of magnetic skyrmion based on the {SAV} approach.
\newblock Submitted.

\bibitem{JiangMSZhangZY22:110954}
M.~Jiang, Z.~Zhang, and J.~Zhao.
\newblock Improving the accuracy and consistency of the scalar auxiliary
  variable ({SAV}) method with relaxation.
\newblock {\em J. Comput. Phys.}, 456:110954, 2022.

\bibitem{JuLLLiX22:66}
L.~Ju, X.~Li, and Z.~Qiao.
\newblock Stabilized exponential-{SAV} schemes preserving energy dissipation
  law and maximum bound principle for the {Allen}–{Cahn} type equations.
\newblock {\em J. Sci. Comput.}, 92(2):66, 2022.

\bibitem{LiPYangL23:182}
P.~Li, L.~Yang, J.~Lan, R.~Du, and J.~Chen.
\newblock A second-order semi-implicit method for the inertial
  {Landau}-{Lifshitz}-{Gilbert} equation.
\newblock {\em Numer. Math. Theor. Meth. Appl.}, 16(1):182--203, 2023.

\bibitem{LiQMeiLQ21:107290}
Q.~Li and L.~Mei.
\newblock Efficient, decoupled, and second-order unconditionally energy stable
  numerical schemes for the coupled {Cahn}–{Hilliard} system in
  copolymer/homopolymer mixtures.
\newblock {\em Comput. Phys. Commun.}, 260:107290, 2021.

\bibitem{LiXShenJ20:2465}
X.~Li and J.~Shen.
\newblock Error analysis of the {Sav}-{Mac} scheme for the {Navier}-{Stokes}
  equations.
\newblock {\em SIAM J. Numer. Anal.}, 58(5):2465--2491, 2020.

\bibitem{LiXLShenJ22:141}
X.~Li, J.~Shen, and Z.~Liu.
\newblock New {SAV}-pressure correction methods for the {Navier}-{Stokes}
  equations: stability and error analysis.
\newblock {\em Math. Comp.}, 91(333):141--167, 2022.

\bibitem{LiXLWangWL22:1026}
X.~Li, W.~Wang, and J.~Shen.
\newblock Stability and error analysis of {IMEX} {SAV} schemes for the
  magneto-hydrodynamic equations.
\newblock {\em SIAM J. Numer. Anal.}, 60(3):1026--1054, 2022.

\bibitem{muMag}
{Micromagnetic Modeling Activity Group, National Institute of Standardsand
  Technology}.
\newblock Available: \url{http://www.ctcms.nist.gov/~rdm/mumag.html}.

\bibitem{OikawaTYokotaH16:056006}
T.~Oikawa, H.~Yokota, T.~Ohkubo, and K.~Hono.
\newblock Large-scale micromagnetic simulation of {Nd-Fe-B} sintered magnets
  with {Dy-rich} shell structures.
\newblock {\em AIP Adv.}, 6:056006, 2016.

\bibitem{RaveWHubertA00:3886}
W.~Rave and A.~Hubert.
\newblock Magnetic ground state of a thin-film element.
\newblock {\em IEEE Trans. Magn.}, 36(6):3886--3899, 2000.

\bibitem{SchafferSSchreflT23:170761}
S.~Schaffer, T.~Schrefl, H.~Oezelt, A.~Kovacs, L.~Breth, N.~J. Mauser,
  D.~Suess, and L.~Exl.
\newblock Physics-informed machine learning and stray field computation with
  application to micromagnetic energy minimization.
\newblock {\em J. Magn. Magn. Mater.}, 576:170761, 2023.

\bibitem{ScholzWFidlerJ03:366}
W.~Scholz, J.~Fidler, T.~Schrefl, D.~Suess, R.~Dittrich, H.~Forster, and
  V.~Tsiantos.
\newblock Scalable parallel micromagnetic solvers for magnetic nanostructures.
\newblock {\em Comput. Mater. Sci.}, 28:366--383, 2003.

\bibitem{ShenJXuJ18:2895}
J.~Shen and J.~Xu.
\newblock Convergence and error analysis for the scalar auxiliary variable
  ({SAV}) schemes to gradient flows.
\newblock {\em SIAM J. Numer. Anal.}, 56(5):2895--2912, 2018.

\bibitem{ShenJXuJ18:407}
J.~Shen, J.~Xu, and J.~Yang.
\newblock The scalar auxiliary variable ({SAV}) approach for gradient flows.
\newblock {\em J. Comput. Phys.}, 353:407--416, 2018.

\bibitem{SussDSchreflT00:3282}
D.~S\"{u}ss, T.~Schrefl, and J.~Fidler.
\newblock Micromagnetics simulation of high energy density permanent magnets.
\newblock {\em IEEE Trans. Magn.}, 36(5):3282--3284, 2000.

\bibitem{TanakaTFuruyaA17:7100304}
T.~Tanaka, A.~Furuya, Y.~Uehara, K.~Shimizu, J.~Fujisaki, T.~Ataka, and
  H.~Oshima.
\newblock Speeding up micromagnetic simulation by energy minimization with
  interpolation of magnetostatic field.
\newblock {\em IEEE Trans. Magn.}, 53(6):7100304, 2017.

\bibitem{TellezMoraAHeX24:1}
A.~Tellez-Mora, X.~He, E.~Bousquet, L.~Wirtz, and A.~H. Romero.
\newblock Systematic determination of a material's magnetic ground state from
  first principles.
\newblock {\em npj Comput. Mater.}, 10(1):1--10, 2024.

\bibitem{WangTZhouJ23:719}
T.~Wang, J.~Zhou, and G.~Hu.
\newblock An {SAV} method for imaginary time gradient flow model in density
  functional theory.
\newblock {\em Adv. Appl. Math. Mech.}, 15(3):719--736, 2023.

\bibitem{YangLChenJR21:110142}
L.~Yang, J.~Chen, and G.~Hu.
\newblock A framework of the finite element solution of the
  {Landau}-{Lifshitz}-{Gilbert} equation on tetrahedral meshes.
\newblock {\em J. Comput. Phys.}, 431:110142, 2021.

\bibitem{YangLHuGH19:1048}
L.~Yang and G.~Hu.
\newblock An adaptive finite element solver fordemagnetization field
  calculation.
\newblock {\em Adv. Appl. Math. Mech.}, 11(5):1048--1063, 2019.

\bibitem{YangXHeX22:114376}
X.~Yang and X.~He.
\newblock A fully-discrete decoupled finite element method for the conserved
  {Allen}–{Cahn} type phase-field model of three-phase fluid flow system.
\newblock {\em Comput. Methods Appl. Mech. Eng.}, 389:114376, 2022.

\bibitem{YaoCDuZ23:139}
C.~Yao, Z.~Du, and L.~Yang.
\newblock {SAV} finite element method for the {Peng}-{Robinson} equation of
  state with dynamic boundary conditions.
\newblock {\em Adv. Appl. Math. Mech.}, 15(1):139--158, 2023.

\bibitem{ZhanJYangL24:1327}
J.~Zhan, L.~Yang, R.~Du, and Z.~Cui.
\newblock Towards preserving geometric properties of
  {Landau}-{Lifshitz}-{Gilbert} equation using multistep methods.
\newblock {\em Commun. Comput. Phys.}, 35:1327--1351, 2024.

\bibitem{ZhuangQQShenJ19:72}
Q.~Zhuang and J.~Shen.
\newblock Efficient {SAV} approach for imaginary time gradient flows with
  applications to one- and multi-component {Bose}-{Einstein} {Condensates}.
\newblock {\em J. Comput. Phys.}, 396:72--88, 2019.

\end{thebibliography}
 \end{document}